\documentclass[letterpaper, 10 pt, conference]{ieeeconf}  

\IEEEoverridecommandlockouts

\title{\LARGE \bf Robust Spacecraft Low-Thrust Trajectory Design: A Chance-Constrained Covariance-Steering Approach}

\author{Meysam Babapour$^{1}$, Ehsan Taheri$^{1}$
\thanks{$^{1}$Meysam Babapour and Ehsan Taheri are with the Department of Aerospace Engineering, Auburn University, Auburn, AL USA,   {\tt\small \{mzb0218,ezt0028\}@auburn.edu}.}%
}

\usepackage{graphicx}
\graphicspath{ {./Figures/}}
\usepackage{amsmath,amsfonts,amssymb,tabularx,xcolor}
\usepackage{bm}
\usepackage{physics}
\usepackage{xcolor}
\usepackage[capitalise]{cleveref}
\usepackage{float}
\usepackage{comment}
\usepackage{subcaption}
\makeatletter
\renewcommand*\env@matrix[1][\arraystretch]{%
  \edef\arraystretch{#1}%
  \hskip -\arraycolsep
  \let\@ifnextchar\new@ifnextchar
  \array{*\c@MaxMatrixCols c}}
\makeatother

\begin{document}

\maketitle
\thispagestyle{empty}
\pagestyle{empty}

\begin{abstract}
This paper proposes a systematic method for generating practical and robust low-thrust spacecraft trajectories. One contribution is to consider the change in mass of the spacecraft at two levels: a) the propulsive acceleration and b) the intensity of the stochastic disturbances. A \textit{covariance variable} formulation is considered, which is computationally more efficient than the factorized covariance implementation. The proposed approach is applied to two- (i.e., planar) and three-dimensional heliocentric phases of spacecraft flight from Earth to Mars under the restricted two-body dynamics. The results highlight the importance of keeping track of mass change to generate more realistic, robust trajectories for interplanetary space missions to avoid underestimation of mission risks.

\end{abstract}

\section{INTRODUCTION}
Optimal control principles (OCPs) are used extensively for spacecraft trajectory design \cite{trelat2012optimal,taheri2018generic,taheri2023l2}. The study of OCPs under uncertainty has a long history. The early developments were grounded in problems solved through dynamic programming (DP), a rigorous framework for sequential decision-making, but severely limited in its applicability to realistic high-dimensional systems \cite{Bellman1966,lee2004approximate}. To address computational barriers, approximate methods such as differential dynamic programming (DDP) were developed, which iteratively refine control laws by local quadratic approximations of the so-called value function \cite{jacobson1970differential}. These methods reduce computational complexity,  but faces difficulties when nonlinear dynamics, nonconvex constraints, or non-Gaussian uncertainties have to be considered \cite{Ozaki2020}.

For unconstrained linear stochastic systems, the state mean and covariance evolve independently. The mean trajectory is governed by the feedforward control inputs, while the covariance dynamics is shaped by the choice of state feedback gains. This decoupling property underlies the observation that solutions to covariance-steering problems with quadratic costs are closely related to Linear Quadratic Regulator (LQR)-based formulations, where feedback laws naturally regulate state dispersion \cite{Chen2016,Chen2016II,Chen2018}.

The independence between mean and covariance no longer holds once constraints are present, and thus mean and covariance dynamics become tightly coupled. As a result, standard LQR methods are not sufficient to address the complexity of practical problems. Alternative approaches have sought to integrate LQR with stochastic optimization. One such strategy consists of designing a feedback gain a priori using LQR and then optimizing only the feedforward control component within a stochastic framework \cite{Lew2020}. A more flexible formulation treats the feedback gain as an optimization variable, enabling joint design of feedforward and feedback terms to better handle nonlinearities, nonconvex and probabilistic constraints, and covariance propagation. The resulting non-convex problems are typically solved using sequential convex programming (SCP) \cite{nurre2022comparison}.

The \textit{factorized covariance} formulation method is proposed in \cite{NonlinearDrag, ridderhof2020chance, DistributionallyRobust} that applies an SCP approach to solve covariance-steering problems. However, factorized covariance formulations are computationally expensive, as each subproblem involves a linear matrix inequality (LMI) whose dimensions grow quadratically with the number of nodes, $N$, leading to a complexity of $\mathcal{O}(N^2 n_x n_u)$ (with $n_x$ and $n_u$ denoting the dimensions of state and control vectors, respectively). An alternative \textit{covariance variable} SCP-based formulation was introduced in \cite{Benedikter2022}. Nevertheless, the proof of lossless convexification underpinning this method was later shown to be incorrect \cite{Rapakoulias2023}. A subsequent modification by \cite{Kumagai2025} (based on the original work \cite{liu2024optimal}) is also  proposed to restore the validity of the formulation under the assumption that a set of tight tolerances are satisfied. Covariance-steering formulations are used for various applications, including path planning \cite{okamoto2019optimal}, powered descent guidance \cite{ridderhof2021minimum}, aerocapture \cite{ridderhof2022chance} and spacecraft rendezvous and docking maneuvers \cite{zhang2023stochastic}.

\textbf{Contributions}: 1) We investigate the application of the \textit{covariance variable} SCP-based formulation for solving fuel-optimal low-thrust Earth-to-Mars problems, which explicitly models spacecraft mass dynamics. This formulation a) allows for thrust to be considered as the control input, b) models a mass-varying stochastic force intensity, and c) allows for enforcing a standard deviation in spacecraft final mass, which ties it nicely with some of the recent research on desensitization of spacecraft trajectories \cite{jawaharlal2024reduced}. In previous applications of covariance steering for low-thrust trajectory design, the mass is not considered as a state. However, the instantaneous mass is an essential component of the propulsive acceleration term that appears in the equations of motion. The advantage of this model is that we can consider realistic engine parameters (such as thrust magnitude and specific impulse) and enforce constraints directly on the maximum thrust that is to be produced by the propulsion system. 2) We also review the steps for implementing a \textit{covariance variable} formulation, mirroring the steps outlined in Ref. \cite{Kumagai2025}, which is computationally more efficient than the \textit{factorized covariance} implementation. However, we provide particular modifications needed for considering the variation in mass. A detailed analysis of the results is also presented.
The paper is structured with Sec. \ref{sec:probform} introducing the general formulation of the stochastic OCP. Section \ref{sec:enforceCOE} describes the detailed convexification procedure and the iterative algorithm. The proposed approach is then assessed in Section \ref{sec:results} through numerical simulations on more realistic scenarios. Finally, Section \ref{sec:con} presents concluding remarks with potential directions for future research.

The notation used in this paper is standard. In particular, $\mathbb{E}[.]$ denotes the expectation operator, $\mathbb{P}(\cdot)$ denotes the probability measure, and $\mathbb{S}^{n}_{+}$ is a set of all positive definite matrices of dimension $n$. Also, $x_t = x(t)$ and $u_t = u(t)$. $\mathbb{R^+}$ denotes the set of positive real numbers. $\bm{I}_{n}$ denotes an $n$-dimensional identity matrix.

\section{Chance-Constrained Stochastic Optimal Control} \label{sec:probform}
\subsection{Continuous-Time Formulation}
We consider a continuous-time nonlinear stochastic dynamical system governed by the stochastic differential equation (SDE), which can be written as,
\begin{equation}
d\boldsymbol{x}_t \;=\; \boldsymbol{f}(\boldsymbol{x}_t,\boldsymbol{u}_t,t)\,dt 
\;+\; \boldsymbol{g}(\boldsymbol{x}_t,\boldsymbol{u}_t)\,d\boldsymbol{w}_t, 
\quad t \in [t_0,t_f],
\label{eq:sde}
\end{equation}
where $\boldsymbol{x}_t \in \mathbb{R}^{n_x}$ denotes the system state, 
$\boldsymbol{u}_t \in \mathbb{R}^{n_u}$ is the control input, and 
$\boldsymbol{f}:\mathbb{R}^{n_x}\times\mathbb{R}^{n_u}\times\mathbb{R}\to \mathbb{R}^{n_x}$ 
represents the drift dynamics. The diffusion term 
$\boldsymbol{g}:\mathbb{R}^{n_x}\times\mathbb{R}^{n_u}\to \mathbb{R}^{n_x\times n_w}$ 
maps the disturbance to the state space, where $\boldsymbol{w}_t \in \mathbb{R}^{n_w}$ 
is a standard Wiener process with independent increments satisfying 
$\mathbb{E}[d\boldsymbol{w}_t]=\boldsymbol{0}_{n_w}$ and 
$\mathbb{E}[d\boldsymbol{w}_t d\boldsymbol{w}_t^\top] = \boldsymbol{I}_{n_w}\,dt$.

Our goal is to steer the trajectory of the system (\ref{eq:sde}) from a normally distributed initial condition $\boldsymbol{x}_{t_0}$ to a desired/specified final distribution $\boldsymbol{x}_{t_f}$
\begin{equation}
\boldsymbol{x}_{t_0}\sim\mathcal{N}(\bar{\boldsymbol{x}}_i, {P}_i), \quad\boldsymbol{x}_{t_f}\sim\mathcal{N}(\bar{\boldsymbol{x}}_f, P_f),
\label{eq:init&finalcondi}
\end{equation}
where $\bar{\boldsymbol{x}}_i,\bar{\boldsymbol{x}}_f \in \mathbb{R}^{n_x}$ are the initial and desired (final)  mean states, and ${P}_i,{P}_f \in \mathbb{S}^{n_x}_{+}$ are the corresponding covariance matrices. To achieve this, we seek a control policy $\boldsymbol{u}_t = \pi(\boldsymbol{x}_t,t)$ that simultaneously regulates both the mean state trajectory and its  covariance evolution. We employ a piecewise constant affine control of the form
\begin{equation}
\boldsymbol{u}_t = \boldsymbol{F}_t +{K}_t\big(\boldsymbol{x}_t - \bar{\boldsymbol{x}}_t\big), \quad t \in [t_0,t_f],
\label{eq:cntrl}
\end{equation}
where $\boldsymbol{F}_t \in \mathbb{R}^{n_u}$ denotes the nominal (feed-forward) control sequence, 
${K}_t \in\mathbb{R}^{n_u \times n_x}$ is the feedback gain matrix.

Throughout the time horizon, the state and control trajectories must satisfy various operational and safety requirements. Since the state distribution and consequently the control input, as a state-feedback law, are unbounded under Gaussian assumptions, these conditions have to be imposed probabilistically as chance constraints, which can be written as,
\begin{equation}
\mathbb{P}\big( \boldsymbol{h}(\boldsymbol{x}_t, \boldsymbol{u}_t) \leq 0 \big) \;\geq\; \beta,
\quad t \in [t_0,t_f],
\label{eq:chancecon}
\end{equation}
where $0 < \beta < 1$ denotes the desired confidence level.

In practice, the control authority is limited, $\|\boldsymbol{u}_t\|\leq u_{\max}$. This constraint is enforced by specializing the general chance constraint \eqref{eq:chancecon} to the control input as,
\begin{equation}
\mathbb{P}\!\left( \|\boldsymbol{u}_t\| \leq u_{\max} \right) \;\geq\; \beta_u, 
\quad t \in [t_0,t_f],
\label{eq:cntrl_cons}
\end{equation}
where $u_{\max}>0$ is the maximum available control effort and $\beta_u$ is its associated confidence level. Please refer to Ref. \cite{okamoto2019input} for incorporating hard constraints into covariance-steering problems. 

While satisfying all constraints, the control objective is also to minimize the cumulative control effort. A Lagrange form of the cost functional is written as, 
\begin{equation}
J = \int_{t_0}^{t_f} \|\boldsymbol{u}_t\|\,dt.
\label{eq:cost}
\end{equation}

The performance index \eqref{eq:cost}, however, is not directly tractable, as it involves integration of a stochastic process. To obtain a numerically tractable objective, we reformulate the cost using the $p$-quantile of the stochastic control norm, which accounts for both the mean control effort and its variability, leading to
\begin{equation}
J_1 = Q\left(\int_{t_0}^{t_f} \|\boldsymbol{u}_t\|\,dt; \;p \right), \qquad 0<p<1,
\label{eq:quantilecost}
\end{equation}
where $Q(\cdot;p)$ denotes the quantile operator. The continuous-time nonlinear chance-constrained stochastic optimal control problem can be stated as Problem 1. 

Problem 1 (Chance-Constrained SOC Optimization)
\begin{subequations}
\label{prob:cc-soc}
\begin{align}
\min_{\boldsymbol{u}_t} \quad & J_1 = Q\!\left(\int_{t_0}^{t_f} \|\boldsymbol{u}_t\| \, dt ; \;p \right) \label{eq:costP1}\\
\text{s.t.} \quad & d\boldsymbol{x}_t = \boldsymbol{f}(\boldsymbol{x}_t,\boldsymbol{u}_t,t)\,dt 
    + \boldsymbol{g}(\boldsymbol{x}_t,\boldsymbol{u}_t)\,d\boldsymbol{w}_t, \label{eq:sdeP1}\\
& \boldsymbol{x}_{t_0} \sim \mathcal{N}(\bar{\boldsymbol{x}}_i, \boldsymbol{P}_i), \label{eq:initP1}\\
& \mathbb{P}\!\left(\|\boldsymbol{u}_t\|\leq u_{\max}\right) \geq \beta_u, 
    \quad t\in[t_0,t_f], \label{eq:cc-cntrlP1}\\
& \boldsymbol{x}_{t_f} \sim \mathcal{N}(\bar{\boldsymbol{x}}_f,\boldsymbol{P}_f), \label{eq:finalcondiP1}\\
& \boldsymbol{u}_t = \boldsymbol{F}_t +{K}_t(\boldsymbol{x}_t - \bar{\boldsymbol{x}}_t), \quad t \in [t_0,t_f]. \label{eq:cntrlP1}
\end{align}
\end{subequations}

\subsection{Discrete-Time Deterministic Formulation}
To formulate a tractable optimization problem, the nonlinear stochastic dynamics (\ref{eq:sde}) are approximated by a linear continuous-time, time-varying stochastic system, which is valid under the assumption of small deviations. Linearization around a reference trajectory $(\hat{\boldsymbol{x}}_t, \hat{\boldsymbol{u}}_t)$ gives  
\begin{equation} \label{eq:linSDE}
d\boldsymbol{x}_t = ({A}(t)\boldsymbol{x}_t + {B}(t)\boldsymbol{u}_t + \boldsymbol{c}(t)\big)\,dt + {G}(t)\,d\boldsymbol{w}_t
\end{equation}
where $G(t) = \boldsymbol{g}(\hat{\boldsymbol{x}}_t,\hat{\boldsymbol{u}}_t)$ and the other matrices are   
\begin{align} \label{eq:AandBofLin}
A(t) & = \left.\frac{\partial \boldsymbol{f}}{\partial \boldsymbol{x}}\right|_{(\hat{\boldsymbol{x}}_t,\hat{\boldsymbol{u}}_t)},
& B(t) = \left.\frac{\partial \boldsymbol{f}}{\partial \boldsymbol{u}}\right|_{(\hat{\boldsymbol{x}}_t,\hat{\boldsymbol{u}}_t)},
\end{align}

\begin{equation} \label{eq:cofLin}
\boldsymbol{c}(t) = \boldsymbol{f}(\hat{\boldsymbol{x}}_t,\hat{\boldsymbol{u}}_t,t) - A(t)\hat{\boldsymbol{x}}_t - B(t)\hat{\boldsymbol{u}}_t.
\end{equation}

With the zero-order hold assumption on the control input, the system is then discretized over the time horizon $[t_0,t_f]$, which is partitioned into $N$ uniform segments. The resulting discrete-time, time-varying linear stochastic dynamics are
\begin{equation} \label{eq:discdyn}
\boldsymbol{x}_{k+1} = A_k \boldsymbol{x}_k + B_k \boldsymbol{u}_k + \boldsymbol{c}_k + G_k \boldsymbol{w}_k,
\end{equation}
where $k = 0, 1, ..., N-1$ and the system matrices are obtained from the continuous-time dynamics using the state transition matrix, $\Phi(t, \tau)$, associated with $A(t)$, as follows
\begin{equation} \label{eq:AfromSTM}
A_k = \Phi(t_{k+1},t_k),
\end{equation}
\begin{equation} \label{eq:BkfromSTM}
B_k = \int_{t_k}^{t_{k+1}} \Phi(t_{k+1},\tau)\,B(\tau)\,d\tau,
\end{equation}
\begin{equation} \label{eq:ckfromSTM}
\boldsymbol{c}_k = \int_{t_k}^{t_{k+1}} \Phi(t_{k+1},\tau)\,\boldsymbol{c}(\tau)\,d\tau,
\end{equation}
\begin{equation} \label{eq:QkfromSTM}
Q_k = \int_{t_k}^{t_{k+1}} \Phi(t_{k+1},\tau)\,G(\tau)G(\tau)^{\top}\Phi(t_{k+1},\tau)^{\top}\,d\tau,
\end{equation}
where $\Phi(t, \tau)$ is the solution to the matrix ordinary differential equation \cite{schaub2003analytical} as,
\begin{equation}\label{eq:stmode}
\frac{d}{dt}\Phi(t, \tau) = A(t)\Phi(t, \tau), \quad \Phi(\tau, \tau) = \boldsymbol{I}_{n_x},
\end{equation}
and $\boldsymbol{w}_k\sim\mathcal{N}(\boldsymbol{0}, \boldsymbol{I}_{n_w})$ denotes a sequence of independent and identically distributed (i.i.d.) Gaussian disturbances. The matrix $G_k$ is chosen such that $G_k G_k^{\top} = Q_k$.
The corresponding mean and covariance propagation dynamics of (\ref{eq:discdyn}) are
\begin{equation} \label{eq:meandyn}
\bar{\boldsymbol{x}}_{k+1} = A_k \bar{\boldsymbol{x}}_k + B_k \boldsymbol{F}_k + \boldsymbol{c}_k,
\end{equation}
\begin{equation} \label{eq:covdyn}
{P}_{k+1} = (A_k + B_k K_k)\,{P}_k\,(A_k + B_k K_k)^{\top} + Q_k,
\end{equation}
with the corresponding initial and final conditions from (\ref{eq:init&finalcondi}),
\begin{align}
\bar{\boldsymbol{x}}_0 & = \bar{\boldsymbol{x}}_i,&{P}_0 &= {P}_i,& \bar{\boldsymbol{x}}_N &= \bar{\boldsymbol{x}}_f, & {P}_N = {P}_f.
\end{align}


The probabilistic control constraint \eqref{eq:cntrl_cons} can be reformulated as a conservative deterministic constraint by applying the triangle inequality to the affine control input as,
\begin{equation} \label{eq:detercntrlcons}
\|\boldsymbol{F}_k\| + \sqrt{Q_{\chi^2_{n_u}}(\beta_u)}\;\sqrt{\lambda_{\max}({K}_k {P}_k {K}_k^\top)} \;\leq\; u_{\max},
\end{equation}
where $Q_{\chi^2_{n_u}}(\beta_u)$ is $\beta_u$-quantile of the chi-squared distribution with $n_u$ degrees of freedom, $\lambda_{\max}(\cdot)$ denotes the largest eigenvalue, and ${K_k}{P_k}{K_k}^\top$ is the control covariance matrix.

Similarly, the objective function (\ref{eq:quantilecost}) can be written in a discrete-time deterministic form as,
\begin{equation} \label{eq:discdetercost}
J_2 = \sum_{k=0}^{N-1}\Big(\|\boldsymbol{F}_k\| + \sqrt{Q_{\chi^2_{n_u}}(p)}\;\sqrt{\lambda_{\max}({K}_k {P}_k {K}_k^\top)}\Big).
\end{equation}

Problem 2 (Discrete-Time Chance-Constrained Deterministic Covariance-Steering Optimization).

Let $\bm{z}_2= {\{\bar{\boldsymbol{x}}_k,\boldsymbol{F}_k,{P}_k,{K}_k}\}$, the problem is
\begin{subequations}
\label{prob:cccs_dt}
\begin{align}
\min_{\bm{z}_2} \quad & J_{2} ,\\
\text{s.t.} \quad & 
\bar{\boldsymbol{x}}_{k+1} = A_k \bar{\boldsymbol{x}}_k + B_k \boldsymbol{F}_k + \boldsymbol{c}_k,\\
& {P}_{k+1} = (A_k + B_k K_k)\,{P}_k\,(A_k + B_k K_k)^{\top} + Q_k,\\
& \|\boldsymbol{F}_k\| + \sqrt{Q_{\chi^2_{n_u}}(\beta_u)}\sqrt{\lambda_{\max}({K}_k {P}_k {K}_k^\top)} \leq u_{\max},\\
& \bar{\boldsymbol{x}}_0 = \bar{\boldsymbol{x}}_i, \quad {P}_0 = {P}_i,\\
& \bar{\boldsymbol{x}}_N = \bar{\boldsymbol{x}}_f, \quad  {P}_N = {P}_f,\\
& \boldsymbol{u}_k = \boldsymbol{F}_k +{K}_k(\boldsymbol{x}_k - \bar{\boldsymbol{x}}_k).
\end{align}
\end{subequations}

\section{Sequential Convex Programming} \label{sec:enforceCOE}
We now present an SCP framework to compute a solution to Problem 2 through linearization and convexification.

\subsection{Covariance Propagation}
The covariance propagation dynamics (\ref{eq:covdyn}) is bilinear in ${K}_k$ and ${P}_k$. To convexify it, two auxiliary variables ${U}_k = {K}_k {P}_k \in \mathbb{R}^{n_u \times n_x}$
and ${Y}_k = {K}_k {P}_k {K}_k^\top \in \mathbb{S}^{n_u}_{+}$ are introduced \cite{liu2024optimal}. Substituting these definitions into (\ref{eq:covdyn}) yields the affine recursion as,
\begin{align}
{P}_{k+1} & = A_k {P}_k A_k^\top + A_k {U}_k^\top B_k^\top \nonumber \\ & + B_k {U}_k A_k^\top + B_k {Y}_k B_k^\top + Q_k.
\end{align}

To ensure consistency between ${P}_k$, ${U}_k$, and ${Y}_k$, the Schur-complement LMI constraint is imposed as,
\begin{equation}
\begin{bmatrix} {P}_k & {U}_k^\top \\ {U}_k & {Y}_k \end{bmatrix} \succeq 0.
\end{equation}

\subsection{Formulating Chance-Constraint Control }
The deterministic equation of the probabilistic control constraint (\ref{eq:detercntrlcons}) is nonconvex due to the square root of the spectral term $\sqrt{\lambda_{\max}({Y}_k)}$. To convexify, we first introduce another auxiliary variable ${\tau}_k \in \mathbb{R}^{+}$ as an upper bound of the maximum standard deviation of the control covariance
\begin{equation}
\lambda_{\max}({Y}_k) \;\leq\; {\tau}_k^2,
\end{equation}
and then linearize the quadratic term around a reference value of ${\hat{\tau}}_k$, 
\begin{equation} \label{eq:taulin}
\lambda_{\max}({Y}_k) - {\hat{\tau}}_k^2 - 2{\hat{\tau}}_k({\tau}_k - {\hat{\tau}}_k) \;\leq\; {\zeta}_k,
\end{equation}
where ${\zeta}_k \in \mathbb{R}^{+}$ is an additional slack variable to allow for small discrepancies due to linearization. The resulting convex second-order cone inequality constraint becomes  
\begin{equation}
\|\boldsymbol{F}_k\| + \sqrt{Q_{\chi^2_{n_u}}(\beta_u)}\, \; {\tau}_k \;\leq\; u_{\max}.
\end{equation}

\subsection{Objective Function}  
In addition to replacing $\lambda_{\max}({Y}_k)$ with its upper bound $\tau_k$ in (\ref{eq:discdetercost}), lossless convexification requires that the gradient with respect to ${Y}_k$ be strictly positive \cite{Pilipovsky2024}, i.e., $\frac{\partial J}{\partial{Y}_k} > 0$. This requirement is not satisfied in the current formulation. To guarantee convexity, we augment the cost with a small regularization term involving the trace of ${Y}_k$ \cite{Kumagai2025} as,
\begin{equation}
J_3 = \sum_{k=0}^{N-1}\Big(\|\boldsymbol{F}_k\| + \sqrt{Q_{\chi^2_{n_u}}(p)}\, \; {\tau}_k + \varepsilon_Y\,\mathrm{tr}({Y}_k)\Big).
\end{equation}

Furthermore, the slack variable ${\zeta}_k$ introduced in (\ref{eq:taulin}) must be penalized to minimize the linearization error,
\begin{equation} \label{Jpen}
J_{\mathrm{pen}} = \sum_{k=0}^{N-1}\Big({\zeta}_k + \tfrac{w}{2} {\zeta}_k^2 + \sqrt{w}\,\|{\zeta_k}\|_1\big),
\end{equation}
where $w$ is a penalty coefficient.

Problem 3 (Convex SDP Covariance-Steering Problem).
Let $\bm{z}_3= {\{\bar{\boldsymbol{x}}_k,\boldsymbol{F}_k,{P}_k,{U}_k,{Y}_k,{\tau}_k,{\zeta}_k\}}$, the convex SDP problem (solved within each SCP iteration) is written as,
\begin{subequations}
\label{prob:cccs_convexSDP}
\begin{align}
\min_{\bm{z}_3} \quad & J_{\mathrm{aug}} \;= \; J_3 + J_{\mathrm{pen}},\\
\text{s.t.} \quad & \bar{\boldsymbol{x}}_{k+1} = A_k \bar{\boldsymbol{x}}_k + B_k \boldsymbol{F}_k + \boldsymbol{c}_k, \\
& {P}_{k+1} = A_k {P}_k A_k^\top + A_k {U}_k^\top B_k^\top  \nonumber \\
&  \quad \quad \quad + B_k {U}_k A_k^\top + B_k {Y}_k B_k^\top + Q_k, \\
& \begin{bmatrix} {P}_k & {U}_k^\top \\ {U}_k & {Y}_k \end{bmatrix} \succeq 0,\\
& \|\boldsymbol{F}_k\| + \sqrt{Q_{\chi^2_{n_u}}(\beta_u)}\, \; \tau_k \;\leq \; u_{\max},\\
& \lambda_{\max}({Y}_k) - \hat{\tau}_k^2 - 2\hat{\tau}_k(\tau_k - \hat{\tau}_k) \;\leq\; \zeta_k,\\
& \bar{\boldsymbol{x}}_0 = \bar{\boldsymbol{x}}_i,\;{P}_0 = {P}_i,\;\bar{\boldsymbol{x}}_N = \bar{\boldsymbol{x}}_f,\;{P}_N = {P}_f.
\end{align}
\end{subequations}
    
This convex reformulation removes bilinearities while preserving the structural coupling between covariance and feedback terms. Regularization of ${Y}_k$ in the cost discourages artificial inflation of these terms and ensures conservative yet effective control of uncertainty.  
The procedure for solving Problem 3 is summarized in Algorithm 1 below. \\
Algorithm 1 (Sequential Convex Programming).
\begin{enumerate}
\item Initialization: choose a nominal feasible trajectory $\{\bar{\boldsymbol{x}}_k, \boldsymbol{F}_k,\hat{\tau}_k\}$ and select penalty weights $\varepsilon_Y, w > 0$.
\item Convexification: Linearize and discretize to obtain $\{{A}_k,B_k,\bm{c}_k,Q_k\}$.
\item Solve Convex SDP: solve problem 3 to obtain $\bm{z} = \{\bar{\boldsymbol{x}}_k^*,\boldsymbol{F}_k^*,{P}_k^*, {U}_k^*, {Y}_k^*, \tau_k^*, \zeta_k^*\}$.
\item Reference Update: set $\bm{z}_k \leftarrow \bm{z}_k^*$.
\item Iterate until the solution converges based on $\varepsilon_x$ and $\varepsilon_{\zeta}$.
\end{enumerate}

\section{Numerical Simulations} \label{sec:results}
\subsection{Planar (2D) low-thrust Earth-to-Mars problem}
We consider a planar low-thrust Earth-to-Mars rendezvous problem with zero hyperbolic excess velocities \cite{taheri2016enhanced}. The spacecraft dynamics are modeled in a Sun-centered Cartesian inertial frame restricted to the ecliptic plane written as,
\begin{equation}
\boldsymbol{f}(\boldsymbol{x}_t,\boldsymbol{u}_t,t) =
\begin{bmatrix}
\boldsymbol{v}_t \\[6pt]
-\mu \dfrac{\boldsymbol{r}_t}{\|\boldsymbol{r}_t\|^{3}} + \dfrac{\boldsymbol{u}_t}{m_t} \\[10pt]
-\dfrac{\|\boldsymbol{u}_t\|}{I_\text{sp} g_0}
\end{bmatrix},~\text{with}~ \bm{u}_t=[u_1, u_2]^\top
\end{equation}
\noindent
where $\boldsymbol{x}_t = [\,\boldsymbol{r}_t^\top,\, \boldsymbol{v}_t^\top,\, m_t\,]^\top \in \mathbb{R}^{5}$ is the state vector, with $\boldsymbol{r}_t \in \mathbb{R}^2$ the position, $\boldsymbol{v}_t \in \mathbb{R}^2$ the velocity, and $m_t \in \mathbb{R^+}$ the spacecraft mass. The control input, $\boldsymbol{u}_t \in \mathbb{R}^{2}$, represents the thrust vector. The parameter $\mu$ denotes the Sun’s gravitational parameter, and $I_\text{sp}$ and $g_0$ are the specific impulse and sea-level gravity values, respectively. We consider the inclusion of mass because its change significantly influences both the optimal trajectory and the stochastic behavior of the system, unlike some previous studies \cite{ridderhof2020chance,Benedikter2022,Kumagai2025}. In addition, we can consider the thrust magnitude as the control input, which is more practical as opposed to considering maximum acceleration. We note, however, that formulating and generating spacecraft trajectories based on propulsive acceleration is a reasonable assumption, in particular, for electric low-thrust propulsion systems and over short time horizons \cite{saloglu2024acceleration} or for preliminary studies \cite{arya2023generation}. 
The corresponding SDEs become
\begin{equation}
d
\begin{bmatrix} 
\boldsymbol{r}_t \\[3pt] 
\boldsymbol{v}_t \\[3pt] 
m_t 
\end{bmatrix}
=
\boldsymbol{f}(\boldsymbol{x}_t,\boldsymbol{u}_t,t)~dt
+
\begin{bmatrix}
\boldsymbol{0}_{2 \times 2} \\[6pt]
\dfrac{\gamma}{m_t} I_2 \\[6pt]
\boldsymbol{0}_{1 \times 2}
\end{bmatrix} d\boldsymbol{w}_t,
\end{equation}
where $d\boldsymbol{w}_t \in \mathbb{R}^2$ is a two-dimensional Wiener process with independent increments and $\gamma$ is the force disturbance intensity. The diffusion term scales $\gamma$ by the inverse of the instantaneous spacecraft mass. As $m_t$ decreases, the same level of force uncertainty induces a larger stochastic acceleration disturbance. The continuous-time dynamics are linearized according to Eqs.~\eqref{eq:AandBofLin}-\eqref{eq:cofLin} and subsequently discretized by Eqs.~\eqref{eq:AfromSTM}-\eqref{eq:stmode} over $N = 40$ time intervals. The initial reference trajectory, required for linearization, is obtained by solving a minimum-fuel  deterministic version of the problem using the CasADi solver \cite{andersson2019casadi} following the work in \cite{nurre2025comparison}. The boundary conditions, along with the relevant physical and mission parameters, are summarized in Table~\ref{tab:params}.
\begin{table}[h!]
\centering
\caption{Parameters for the planar Earth-to-Mars problem.}
\begin{tabular}{lcc}
\hline
\textbf{Parameter} & \textbf{Value} & \textbf{Unit} \\
\hline
$\boldsymbol{r}_{i}$ & $[-140699693; -51614428]$ & km \\
$\boldsymbol{v}_{i}$ & $[9.774596; -28.07828]$ & km/s \\
$\boldsymbol{r}_{f}$ & $[-172682023; 176959469]$ & km \\
$\boldsymbol{v}_{f}$ & $[-16.427384; -14.860506]$ & km/s \\
$m_i$ & $5000$ & kg \\
$\sigma_{r_i}$ & $10$ & km \\
$\sigma_{v_i}$ & $0.1$ & km/s \\
$\sigma_{m_i}$ & $0.0$ & kg \\
$\sigma_{r_f}$ & $3.16 \times 10^5$ & km \\
$\sigma_{v_f}$ & $0.1$ & km/s \\
$\sigma_{m_f}$ & $70.7107$ & kg \\
$\mu$ & $1.3271 \times 10^{11}$ & km$^3$/s$^2$ \\
$I_\text{sp}$ & $3000$ & s \\
$g_0$ & $9.80665$ & m/s$^2$ \\
$u_\text{max}$ & $5$ & N \\
$\gamma$ & $9 \times 10^{-5}$ & kg.km/s$^{3/2}$ \\
$t_f$ & $348.795$ & day \\
$d$ & $100$ & -- \\
$\varepsilon_Y$ & $0.01$ & -- \\
$p,\beta_{u}$ & $0.95$ & -- \\
\hline
\end{tabular}
\label{tab:params}
\end{table}
Although the chosen disturbance intensity (i.e., $\gamma/m(t_0)$) is comparable to the magnitude reported in the literature \cite{ridderhof2020chance,Benedikter2022}, the effective stochastic disturbance level exceeds those values by the end of the maneuver due to the decreasing spacecraft mass. Note that $m(t_0)$ is assumed to be deterministic at departure, and a large terminal mass variance is permitted, reflecting the fact that the control objective does not involve regulating the final mass uncertainty (i.e., variance).

To avoid numerical conditioning issues during optimization, all problem data are properly scaled. In particular, the covariance dynamics and related constraints are scaled with a scaling factor of $d$ \cite{Kumagai2025}. The SCP algorithm terminates when the relative difference between consecutive state trajectories satisfies the following constraint as,  
\[
\frac{\|\,\boldsymbol{x}^{(i)} - \boldsymbol{x}^{(i-1)}\,\|}{\|\,\boldsymbol{x}^{(i-1)}\,\|} \leq \varepsilon_x, \quad \varepsilon_x = 5 \times 10^{-4},
\]  
where $(i)$ denotes the SCP iteration and the linearization slack variable satisfies $
\zeta \leq \varepsilon_\zeta = 10^{-6}$. The penalty coefficient $w$ is updated at the $i$-th iteration as,
\[
w = \min\left(10^{(i+3)}, 10^{12}\right).
\]

\begin{figure}[h!]
    \centering
    \begin{subfigure}[b]{0.45\textwidth}
        \centering
        \includegraphics[width=\textwidth]{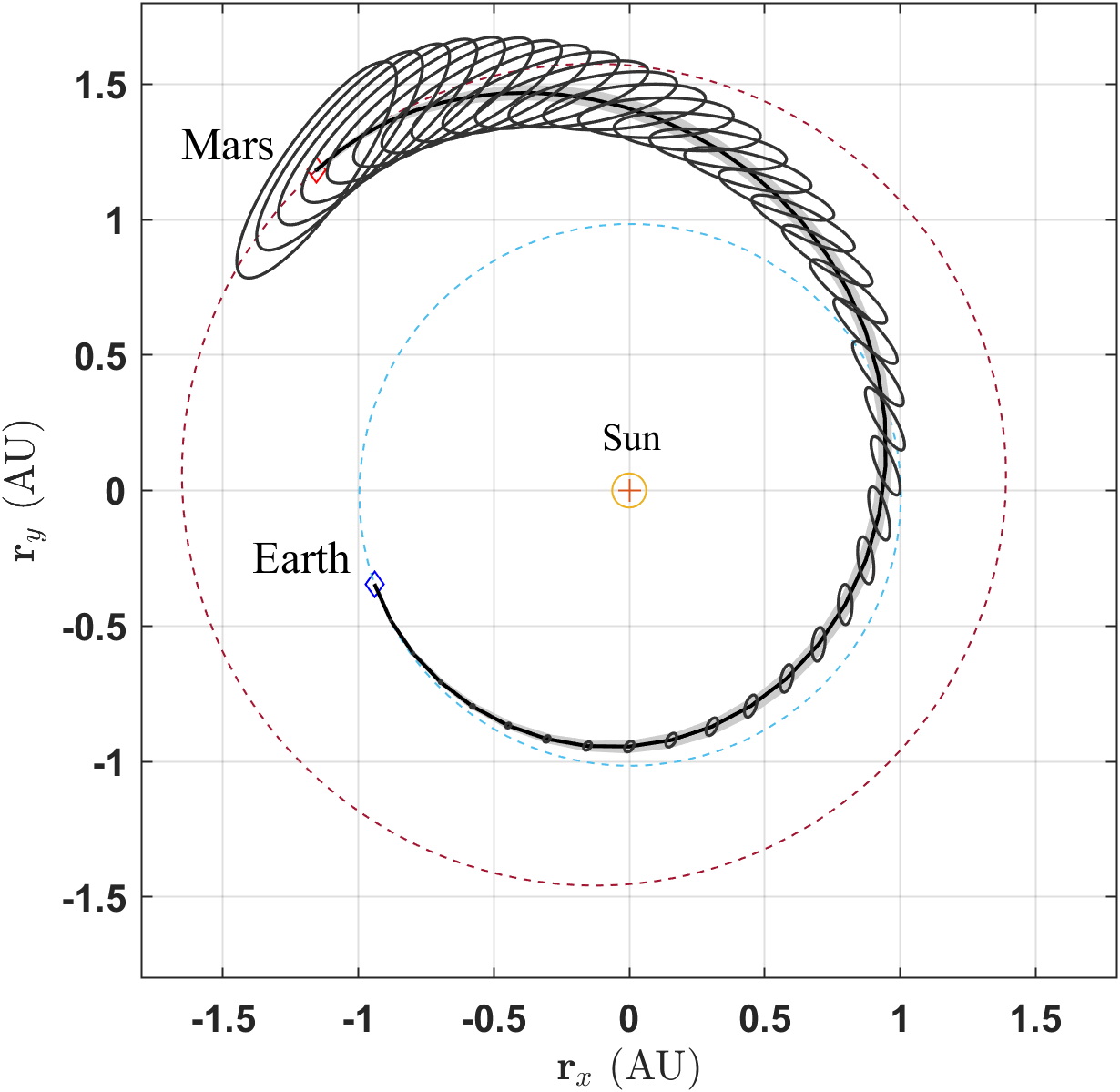}
        \caption{Open-loop position covariance ellipses.}
        \label{fig:OLposcov}
    \end{subfigure}
    \hfill
    \begin{subfigure}[b]{0.45\textwidth}
        \centering
        \includegraphics[width=\textwidth]{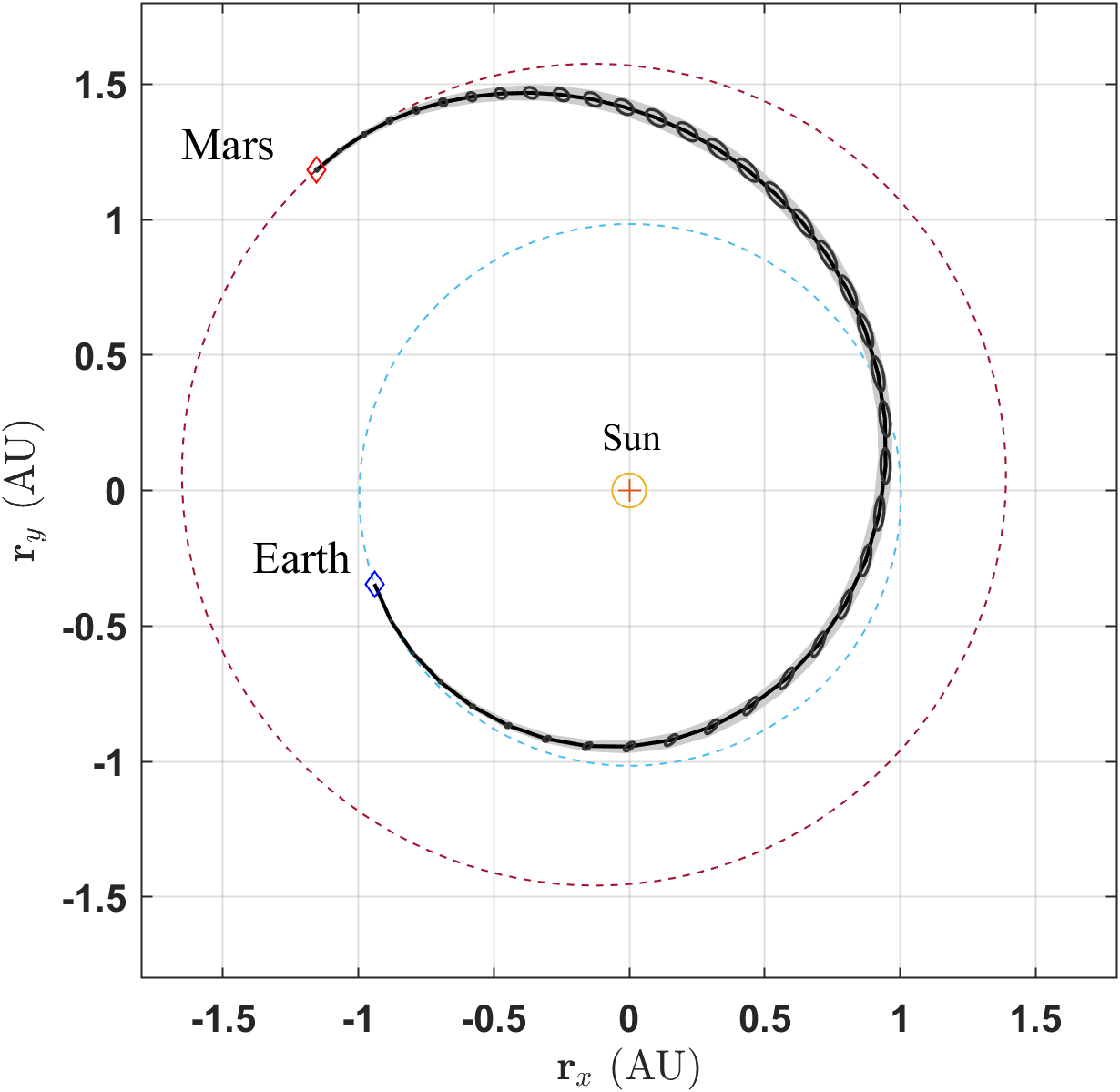}
        \caption{Closed-loop position covariance ellipses.}
        \label{fig:CLposcov}
    \end{subfigure}
    \caption{2D case: open-loop and closed-loop trajectories along with their position covariance ellipses.}
    \label{fig:OL&CLposcov}
\end{figure}

\begin{figure}[h!]
    \centering
    \begin{subfigure}[b]{0.45\textwidth}
        \centering
        \includegraphics[width=\textwidth]{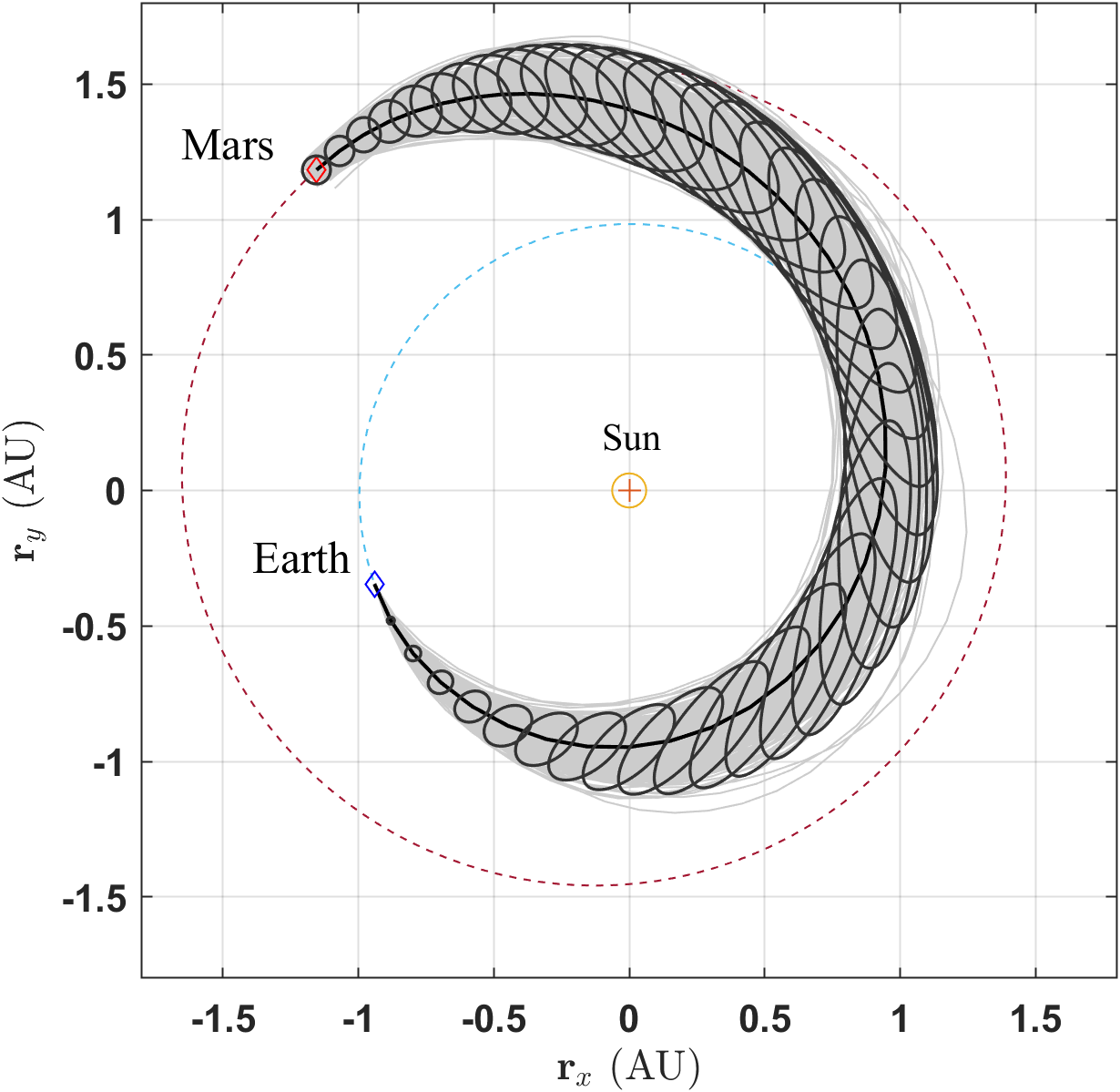}
        \caption{Closed-loop position covariance ellipses.}
        \label{fig:posS10}
    \end{subfigure}
    \hfill
    \begin{subfigure}[b]{0.45\textwidth}
        \centering
        \includegraphics[width=\textwidth]{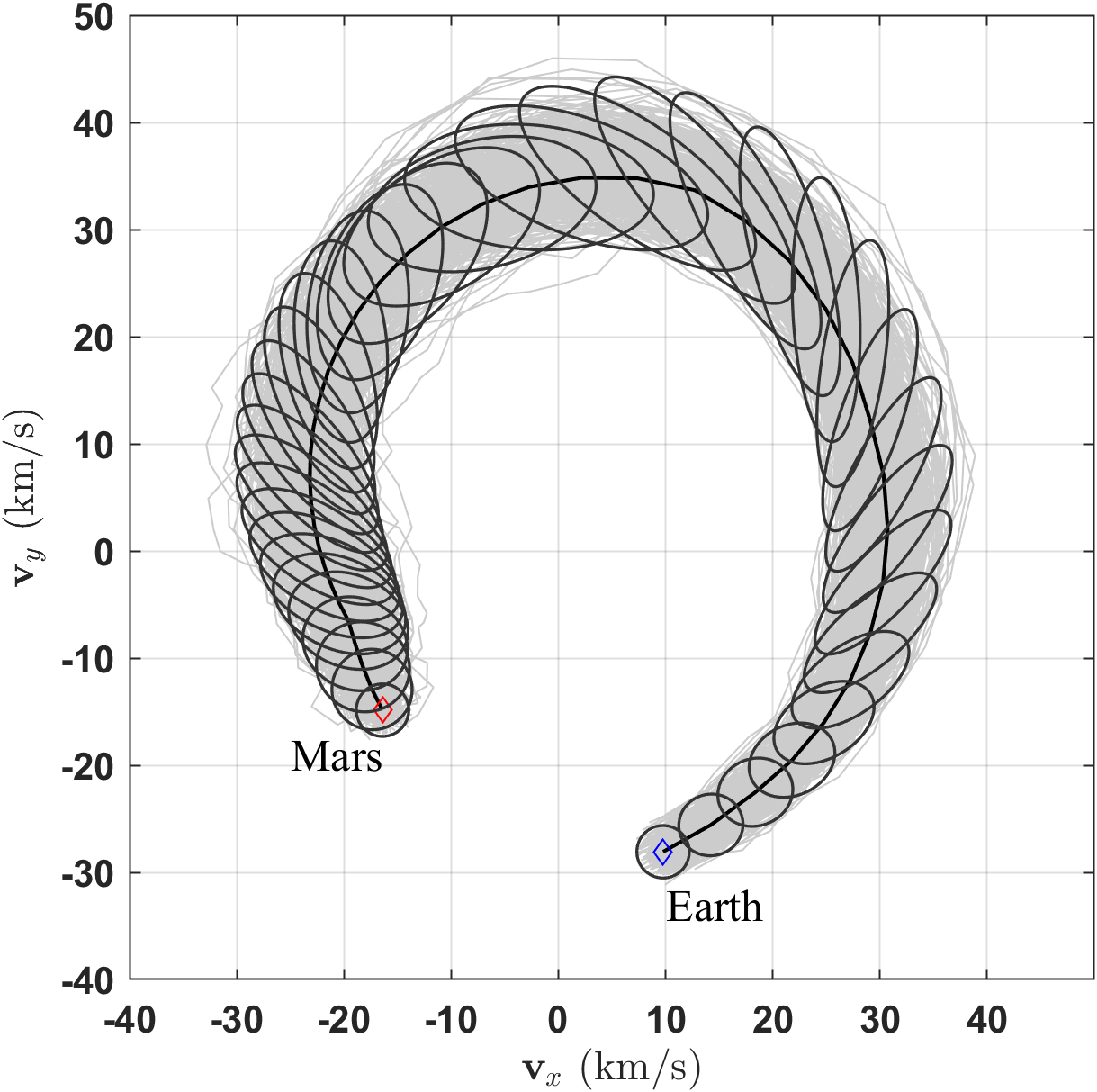}
        \caption{Closed-loop velocity covariance ellipses.}
        \label{fig:velS10}
    \end{subfigure}
    \caption{2D case: position and velocity covariance ellipses (scaled 10x) of the closed-loop solution.}
    \label{fig:CLpos&VelcovS10}
\end{figure}

\begin{figure}[h!]
    \centering
    \begin{subfigure}[b]{0.45\textwidth}
        \centering
        \includegraphics[width=\textwidth]{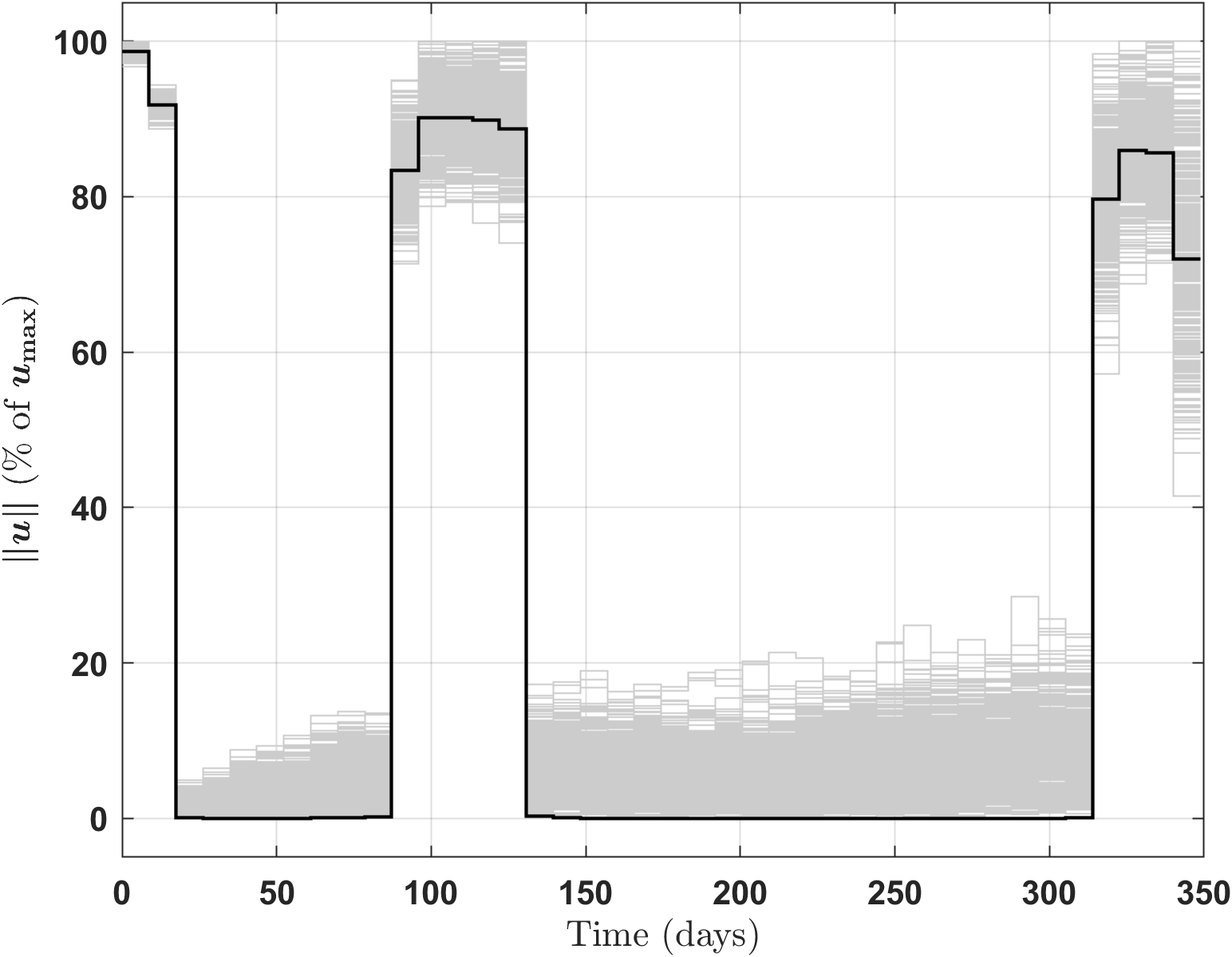}
        \caption{Magnitude of the feedforward control (black) and Monte Carlo control samples (gray).}
        \label{fig:Cntrl}
    \end{subfigure}
    \hfill
    \begin{subfigure}[b]{0.45\textwidth}
        \centering
        \includegraphics[width=\textwidth]{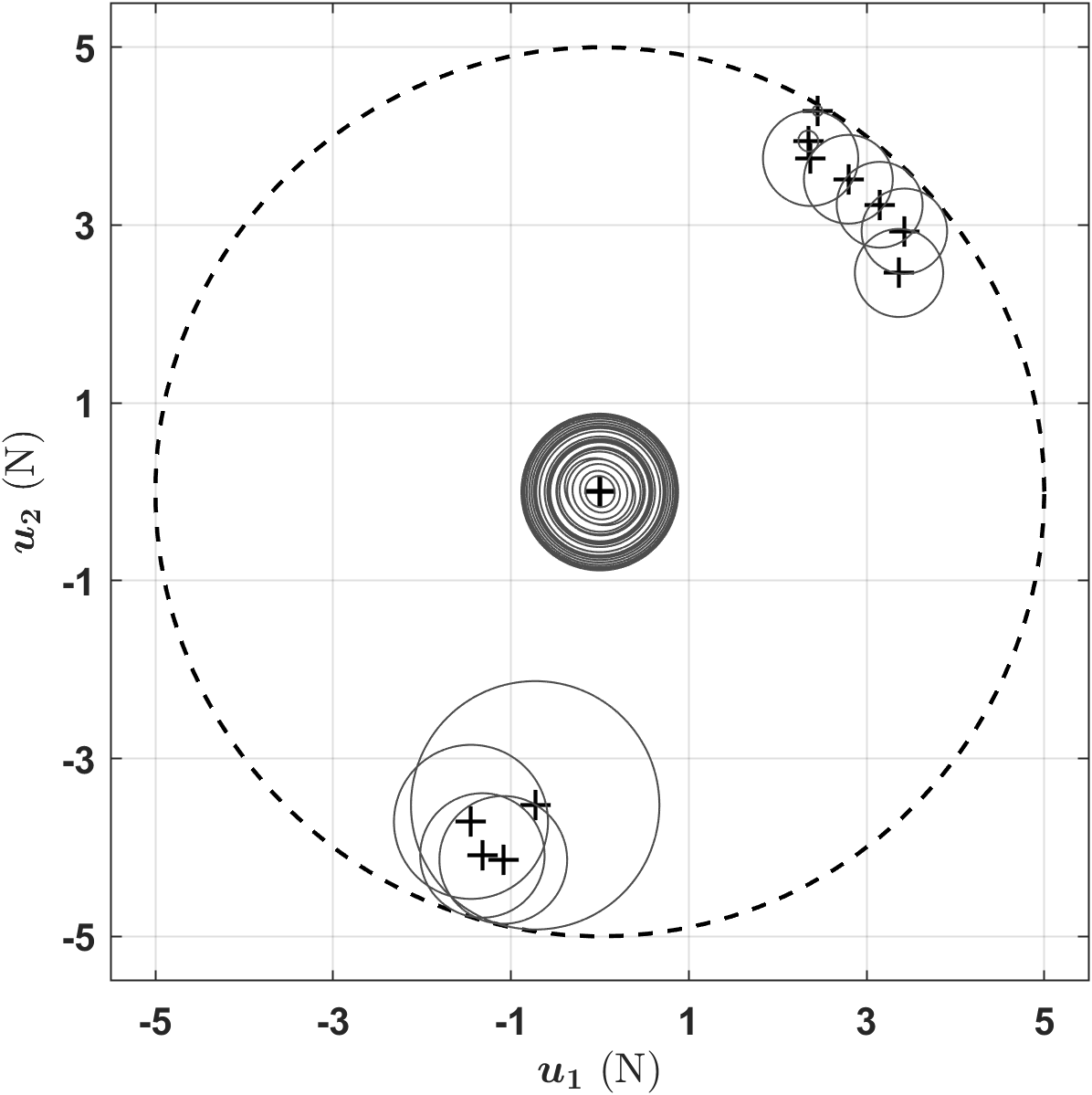}
        \caption{Feedforward control, denoted as `+' with feedback gain uncertainty ellipses.}
        \label{fig:CntrlCov}
    \end{subfigure}
    \caption{2D case: control mean and covariance ellipses.}
    \label{fig:CntrlNom&Cov}
\end{figure}

Problem 3 is solved, and the convergence is achieved in 12 iterations (6.01~s) on an Intel\textregistered~Core\texttrademark~i7-13620H @ 2.40~GHz laptop computer. Figure \ref{fig:OL&CLposcov} depicts the evolution of the position covariance along with the mean trajectory. In the open-loop case, the covariance grows rapidly over time, whereas in the closed-loop case, the $95\%$-confidence ellipses remain bounded due to the action of the feedback controller.
A Monte Carlo simulation (with 1000 samples) was used to validate the optimized trajectory. Figure \ref{fig:CLpos&VelcovS10} compares the resulting dispersion with the nominal solution, in which both position and velocity covariance ellipses in Figs. \ref{fig:posS10} and \ref{fig:velS10} are enlarged tenfold for improved visual interpretation. The close alignment of position and velocity samples within their predicted uncertainty ellipses validates the precision of the covariance propagation framework.
Figure \ref{fig:CntrlNom&Cov} depicts the control profiles and their components. Figure \ref{fig:Cntrl} shows the feed-forward control ($\|\boldsymbol{F} \|$) and Monte Carlo simulation (1000 samples) profiles. While the feedforward control consists of three main thrust arcs (similar to the deterministic solution for the same problem \cite{taheri2016enhanced}), the magnitude of the thrust shows noticeable deviations from a clean bang-bang profile. This is a direct consequence of trading off pure fuel optimality for trajectory robustness, which requires the mean thrust to create some margins for the feedback term to control the covariances, which has occurred mainly during the last two thrust arcs. Figure \ref{fig:CntrlCov} shows the components of the nominal control with feedback gain uncertainty ellipses.

\begin{figure}[h!]
    \centering
    \includegraphics[width = 0.45\textwidth]{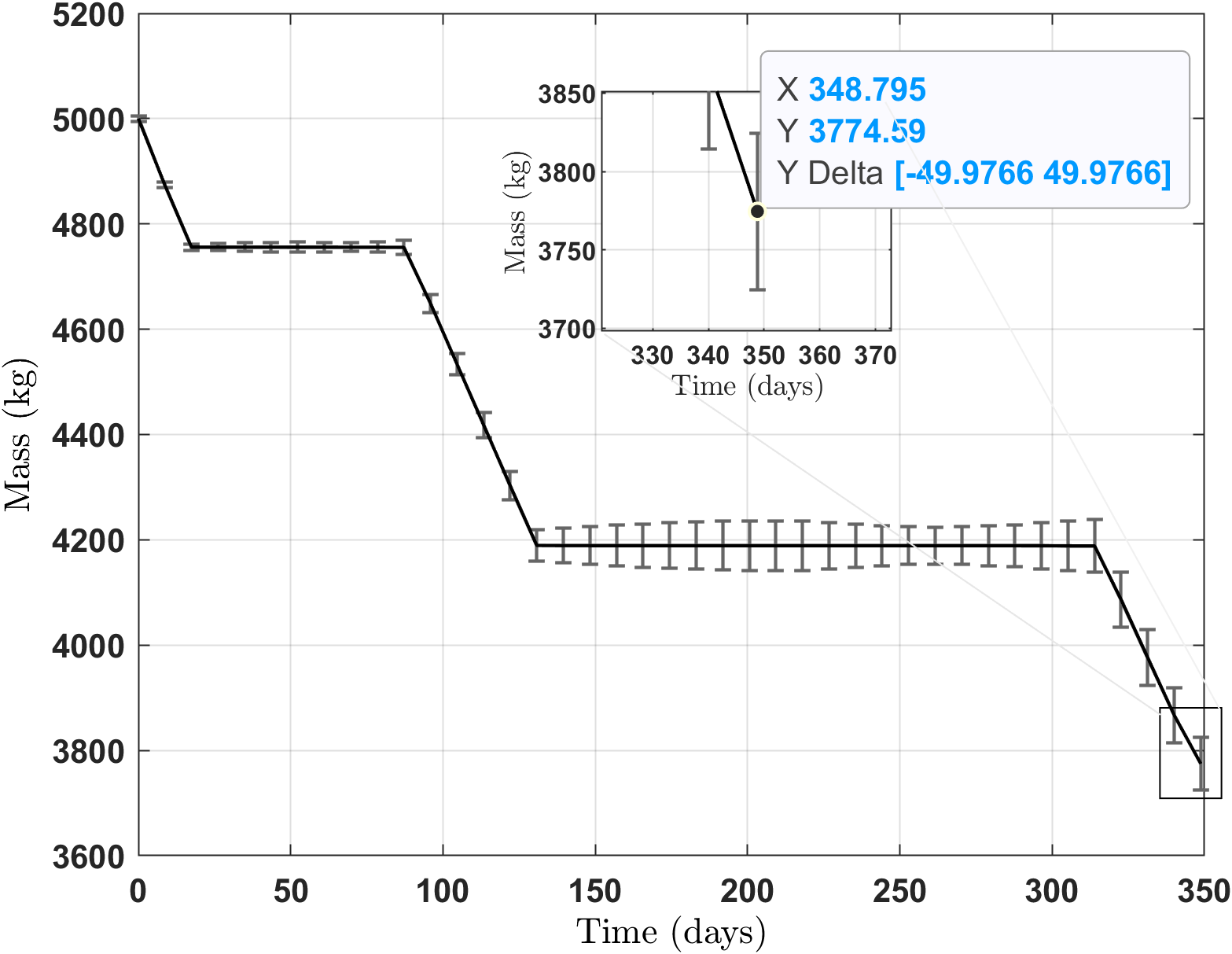}
    \caption{2D case: mass and its standard deviation vs. time.}
    \label{fig:massstd}
\end{figure}

\begin{figure}[h!]
    \centering
    \includegraphics[width = 0.45\textwidth]{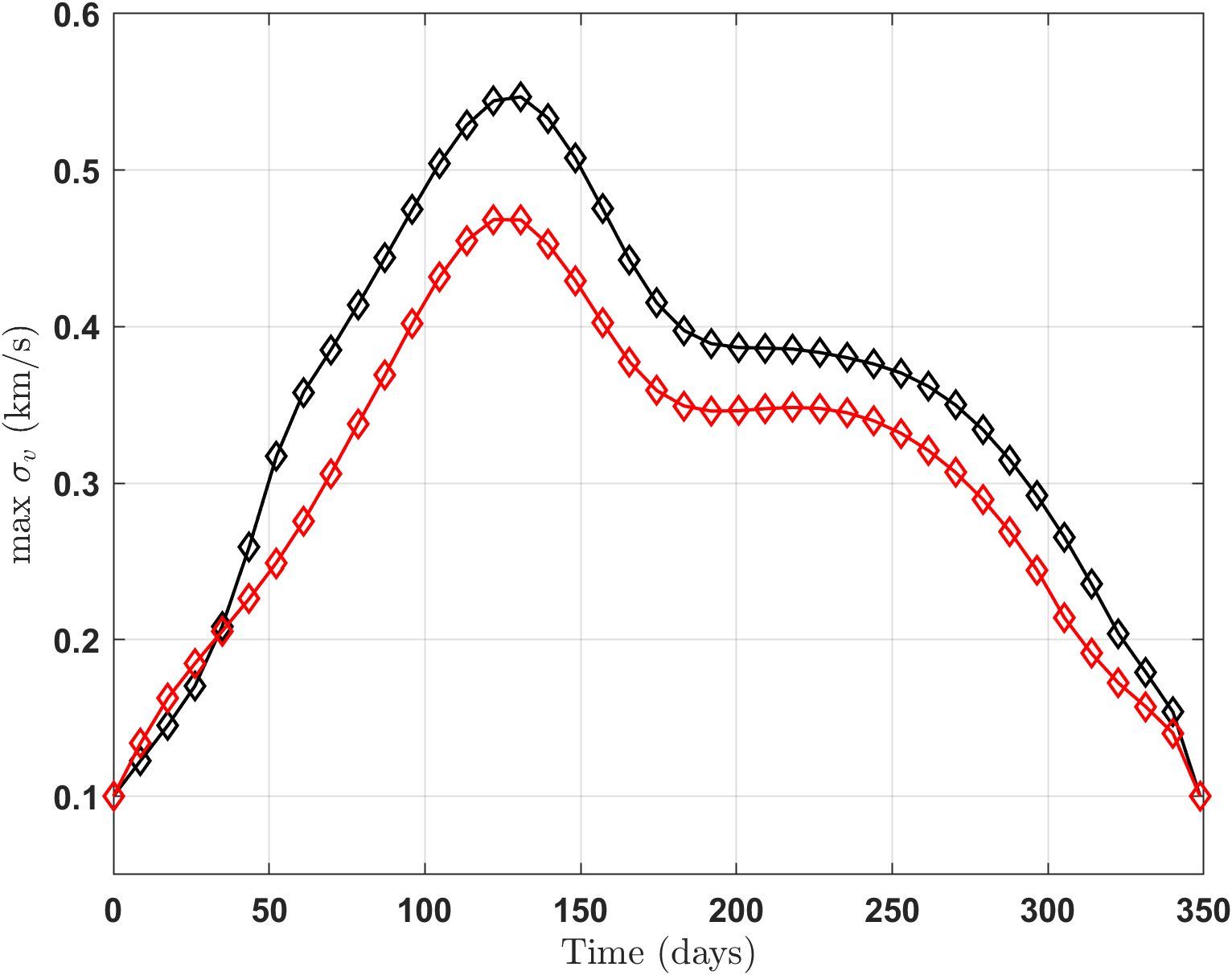}
    \caption{2D case: velocity standard deviation along the first principal component with (black) and without (red) mass variations.}
    \label{fig:velstd}
\end{figure}

\begin{figure}[h!]
    \centering
    \includegraphics[width = 0.45\textwidth]{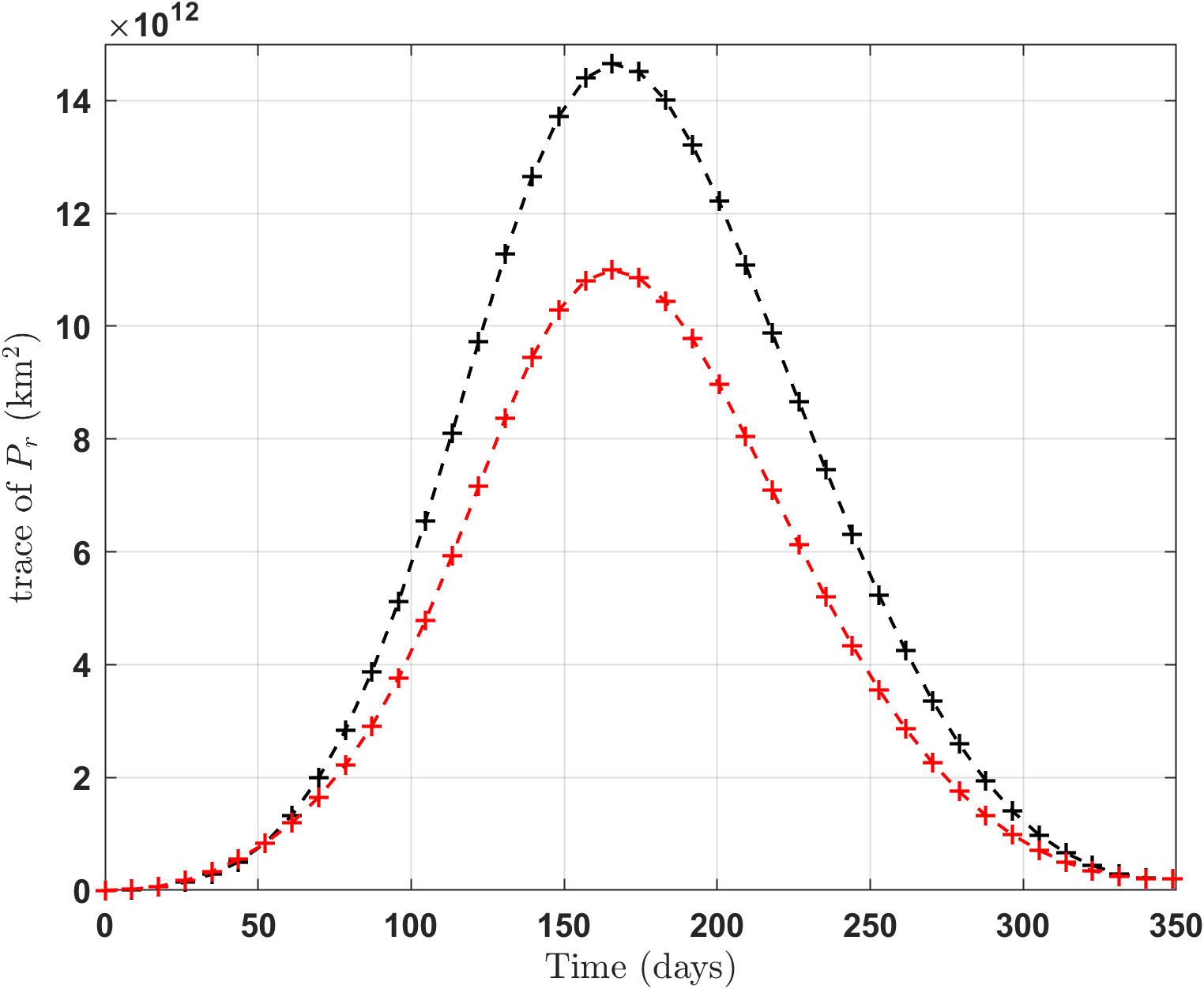}
    \caption{2D case: trace of the position covariances with (black) and without (red) the inclusion of mass variations.}
    \label{fig:postrace}
    \vspace{-4mm}
\end{figure}

Figure \ref{fig:massstd} shows the time histories of the mass and its standard deviation along the trajectory (which is calculated from the covariance solution matrix at the final point). The mass decreases over the course of the mission, reaching a final value of approximately $m_f = 3774.59$ kg upon arrival at Mars. The velocity standard deviation, shown in Figure \ref{fig:velstd}, exhibits a significant growth over the mission when the mass uncertainty is included, with peak standard deviations increasing by up to $16.79\%$. A similar trend is expected in the evolution of the position uncertainty, and as Figure \ref{fig:postrace} shows the trace of the position covariances, the inclusion of mass results in $33.27\%$ higher dispersion at the mid-course phase. These differences underscore the indirect amplification mechanism by which uncertainty in mass, originating from stochastic control input norms, propagates through the dynamics and inflates the velocity variance due to the $1/m$ coupling in the translational acceleration equation. Neglecting this effect leads to an underestimation of the control authority required to robustly satisfy the terminal state covariance constraints, especially in long-duration, low-thrust missions.

\begin{figure}[h]
    \centering
    \includegraphics[width = 0.45\textwidth]{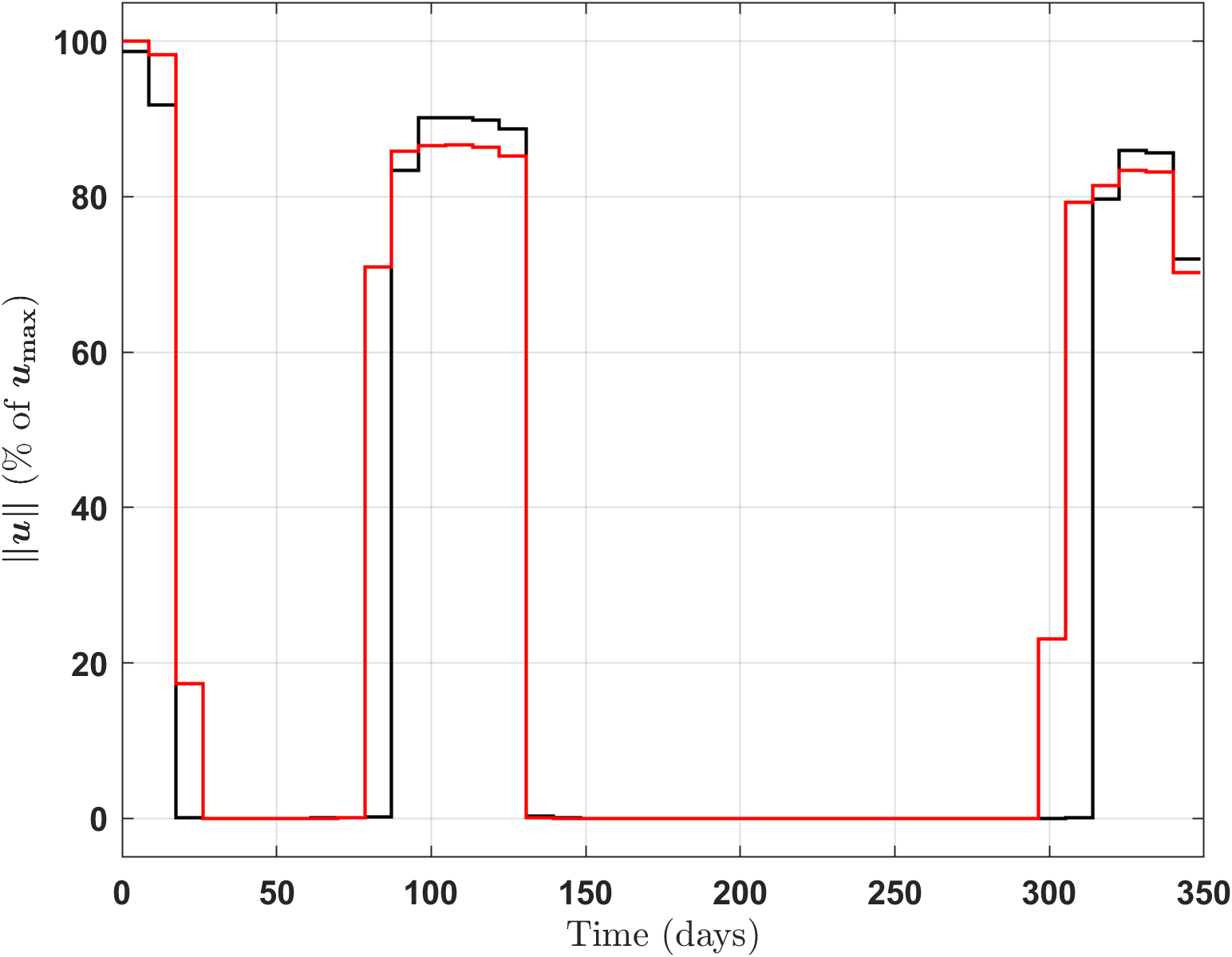}
    \caption{2D case: control w/ (black) and w/o (red) mass uncertainty.}
    \label{fig:twoNominals}
    \vspace{-6mm}
\end{figure}

We highlight that the effect of explicitly modeling mass dynamics, as a stochastic decision variable, was rigorously evaluated by comparing two formulations of the chance-constraint covariance-steering problem. In the first scenario, the mass equation was incorporated as a fifth state, and its change was propagated through the full closed-loop covariance dynamics. In the second scenario, mass evolution was considered deterministic and excluded from the set of decision variables, effectively assuming that the thrust-to-acceleration mapping remained constant throughout the mission. Despite producing comparable nominal trajectories, the results reveal notable discrepancies in the control effort and state uncertainty, particularly during the thrusting phases.
Figure \ref{fig:twoNominals} compares the feedforward control profiles ($\| \bm{F} \|$) for the two scenarios. The results indicate that including mass variations leads to noticeably higher thrust magnitudes during the main thrust arcs, with a peak increase of approximately $6\%$ relative to the case without mass modeling. Note that disturbance becomes stronger in the scenario that mass is decreased over the mission time.

\subsection{Three-dimensional (3D) Earth-to-Mars problem} \label{sec:ch4case3}

This case study presents a more realistic low-thrust Earth-to-Mars rendezvous to demonstrate the full capability of the proposed covariance steering framework.

The spacecraft dynamics are modeled in a Sun-centered Cartesian inertial frame, with the full state vector defined as $\boldsymbol{x}_t = [\,\boldsymbol{r}_t^\top,\, \boldsymbol{v}_t^\top,\, m_t\,]^\top \in \mathbb{R}^{7}$, where $\boldsymbol{r}_t =[r_x,r_y,r_z]^\top \in \mathbb{R}^3$, $\boldsymbol{v}_t =[v_x,v_y,v_z]^\top \in \mathbb{R}^{3}$ are the position and velocity vectors, respectively, and $m_t \in \mathbb{R}^{+}$ is the spacecraft mass.
The dynamics are given by
\begin{equation}
\boldsymbol{f}(\boldsymbol{x}_t,\boldsymbol{u}_t,t) =
\begin{bmatrix}
\boldsymbol{v}_t \\[6pt]
-\mu \dfrac{\boldsymbol{r}_t}{\|\boldsymbol{r}_t\|^{3}} + \dfrac{\boldsymbol{u}_t}{m_t} \\[10pt]
-\dfrac{\|\boldsymbol{u}_t\|}{I_\text{sp} g_0}
\end{bmatrix},
\end{equation}
where $\boldsymbol{u}_t \in \mathbb{R}^{3}$ is the thrust vector, $\mu$ is the Sun's gravitational parameter, and $I_\text{sp}$ and $g_0$ are the specific impulse and standard gravity, respectively.
The mass dynamics are explicitly included to capture the critical coupling between fuel consumption and translational acceleration.

The corresponding stochastic differential equation is
\begin{equation}
d
\begin{bmatrix}
\boldsymbol{r}_t \\[3pt]
\boldsymbol{v}_t \\[3pt]
m_t
\end{bmatrix}
=
\boldsymbol{f}(\boldsymbol{x}_t,\boldsymbol{u}_t,t)~dt
+
\begin{bmatrix}
\boldsymbol{0}_{3 \times 3} \\[6pt]
\frac{\gamma}{m_t} I_3 \\[6pt]
0
\end{bmatrix} d\boldsymbol{w}_t,
\end{equation}
where $d\boldsymbol{w}_t \in \mathbb{R}^{3}$ is a Wiener process and $\gamma$ is the constant acceleration disturbance intensity. The continuous-time dynamics are linearized and discretized with $N = 60$.
The initial nominal trajectory is generated by solving a deterministic minimum-fuel problem using CasADi solver.

The spacecraft is required to move from an initial distribution $\mathcal{N}(\bar{\boldsymbol{x}}_i, P_i)$ to a final distribution $\mathcal{N}(\bar{\boldsymbol{x}}_f, P_f)$.
The initial state uncertainty is set to $\sigma_{r_i} = 10$ km and $\sigma_{v_i} = 0.1$ km/s, with a known initial mass ($\sigma_{m_i} = 0$), resulting in
\begin{align}
P_{i} = \mathrm{diag}(10^2,\ 10^2,\ 10^2,\ 10^{-2}, \nonumber \\ \ 10^{-2},\ 10^{-2},\ 0) \quad \text{(km}^2, \text{(km/s)}^2, \text{kg}^2).
\end{align}

The final state is constrained to $\sigma_{r_f} = 3.16 \times 10^{5}$ km and $\sigma_{v_f} = 0.1$ km/s, with a large value for the variance of the final mass, yielding
\begin{align}
P_{f} = \mathrm{diag}(10^{11},\ 10^{11},\ 10^{11},\ 10^{-2},\ \nonumber \\ 10^{-2},\ 10^{-2},\ 5000) \quad \text{(km}^2, \text{(km/s)}^2, \text{kg}^2).
\end{align}

Please note that a large, unpenalized variance is permitted for the final mass, $\sigma_{m_f} \leq  \sqrt{5000}\approx 70.71$ kg, reflecting that the control objective does not involve regulating the dispersion in the final mass.
The boundary conditions and physical parameters are summarized in Table~\ref{tab:params3D}.

\begin{table}[h!]
\centering
\caption{Parameters for the 3D Earth-to-Mars problem.}
\begin{tabular}{lcc}
\hline
\textbf{Param.} & \textbf{Value} & \textbf{Unit} \\
\hline 
$\bar{\boldsymbol{r}}_{i}$ & $[-140699693; \, -51614428; \, 980]$ & km \\
$\bar{\boldsymbol{v}}_{i}$ & $[9.774596; \, -28.07828; \, 4.337725\times 10^{-4}]$ & km/s \\
$\bar{\boldsymbol{r}}_{f}$ & $[-172682023; \, 176959469; \, 7948912]$ & km \\
$\bar{\boldsymbol{v}}_{f}$ & $[-16.427384; \, -14.860506; \, 9.21486\times 10^{-2}]$ & km/s \\
$\bar{m}_i$ & $5000$ & kg \\
$\mu$ & $1.3271 \times 10^{11}$ & km$^3$/s$^2$ \\
$I_\text{sp}$ & $3000$ & s \\
$g_0$ & $9.80665 \times 10^{-3}$ & km/s$^2$ \\
$u_{\max}$ & $5$ & N \\
$\gamma$ & $9 \times 10^{-5}$ & kg.km/s$^{3/2}$ \\
$t_f$ & $348.795$ & day \\
$\varepsilon_Y$ & $0.01$ & -- \\
$d$ & $100$ & -- \\
\hline
\end{tabular}
\label{tab:params3D}
\end{table}

The SCP termination criteria were set to $\varepsilon_x = 5.0 \times 10^{-4}$ and $\varepsilon_\zeta = 1.0 \times 10^{-6}$. The solver converged to a feasible solution in 12 iterations, demonstrating the method's robustness and scalability to a problem with  7 states, 3 control inputs, and with $N = 60$ grid points.

Unlike the position and velocity states, the initial mass is assumed to be known precisely ($\sigma_{m_i} = 0$), reflecting accurate pre-launch propellant loading measurements.
However, mass becomes a stochastic variable over time because the control input $\boldsymbol{u}_t$, which has an associated execution uncertainty, directly drives the mass rate equation $dm_t/dt \propto -\| \boldsymbol{u}_t \|$. This propagated mass uncertainty then couples back into the translational dynamics.
The state transition matrix (STM) includes the sensitivity of acceleration to mass variations, $\partial (\boldsymbol{u}_t/m_t) / \partial m_t = -\boldsymbol{u}_t/m_t^2$, creating a feedback loop where uncertainty in mass induces uncertainty in the achieved acceleration for a given thrust command and also affects the strength of the disturbance.
\begin{figure}[H]
    \centering
    \includegraphics[width=0.45\textwidth]{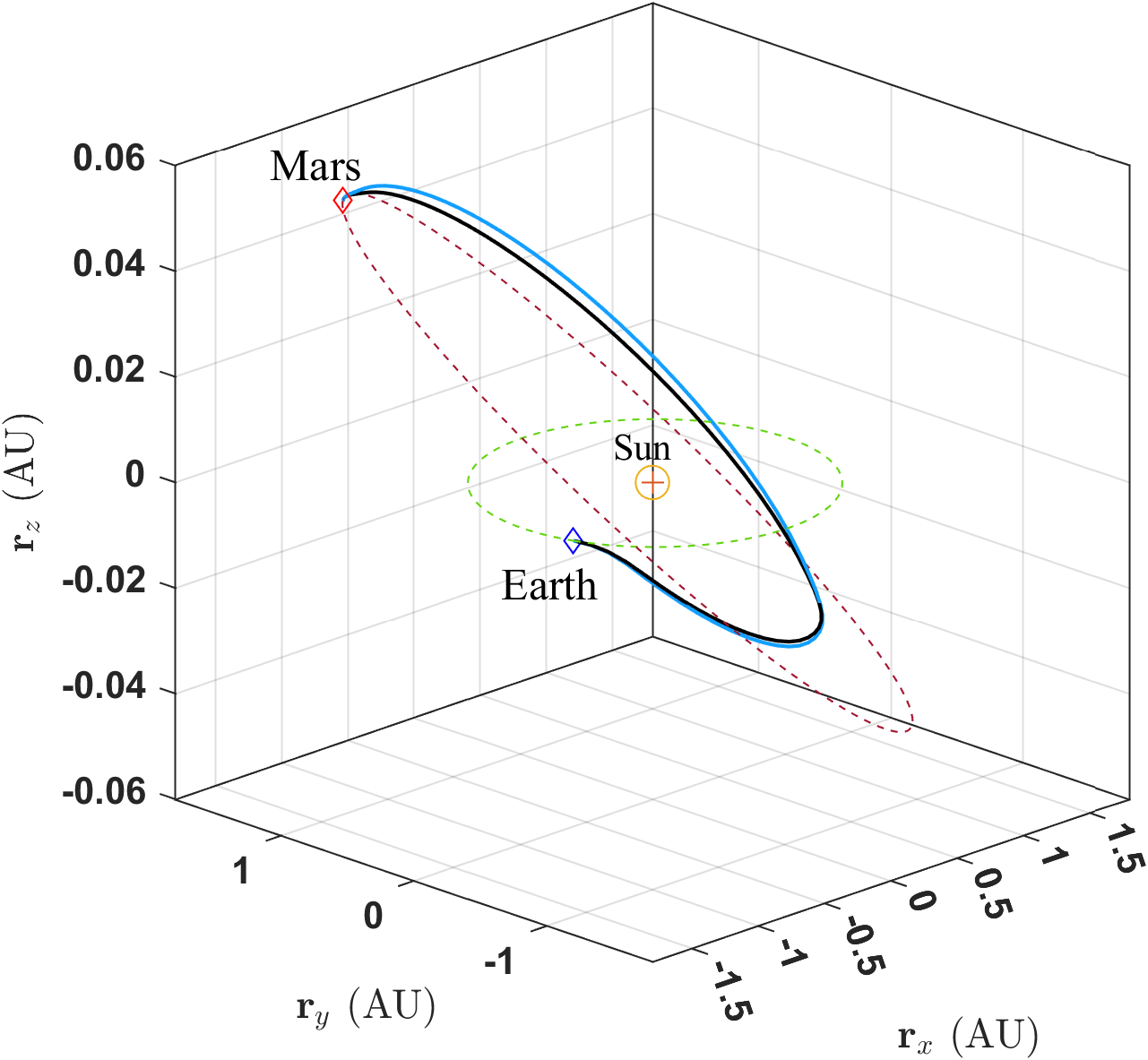}
    \caption{3D view of the deterministic CasADi (black) and closed-loop (blue) trajectories for the Earth-to-Mars problem.}
    \label{fig:3d_traj_cov}
\end{figure}
The optimized controller successfully navigates the uncertain dynamics to meet the stringent terminal constraints. Figure \ref{fig:3d_traj_cov} shows the resulting 3D nominal trajectories, illustrating the spacecraft's path from Earth to Mars in the more realistic 3D dynamical model. A Monte Carlo simulation with 1000 samples validates the controller's performance and the accuracy of the covariance propagation. 

\begin{figure}[H]
    \centering
    \begin{subfigure}[b]{0.5\textwidth}
        \centering
        \includegraphics[width=\textwidth]{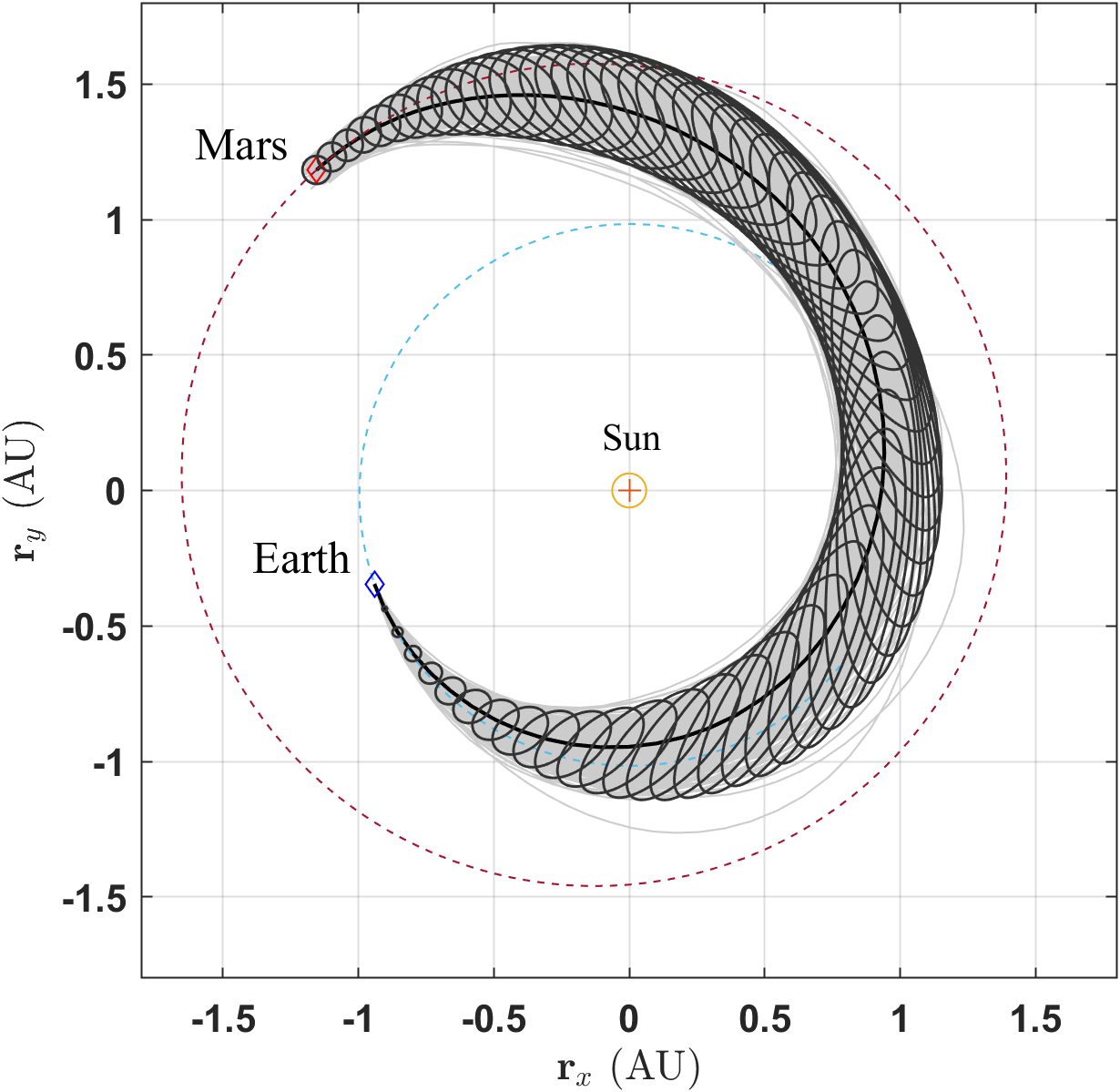}
        \caption{Closed-loop trajectory and position covariances (10x).}
        \label{fig:3DposS10}
    \end{subfigure}
    \hfill
    \begin{subfigure}[b]{0.5\textwidth}
        \centering
        \includegraphics[width=\textwidth]{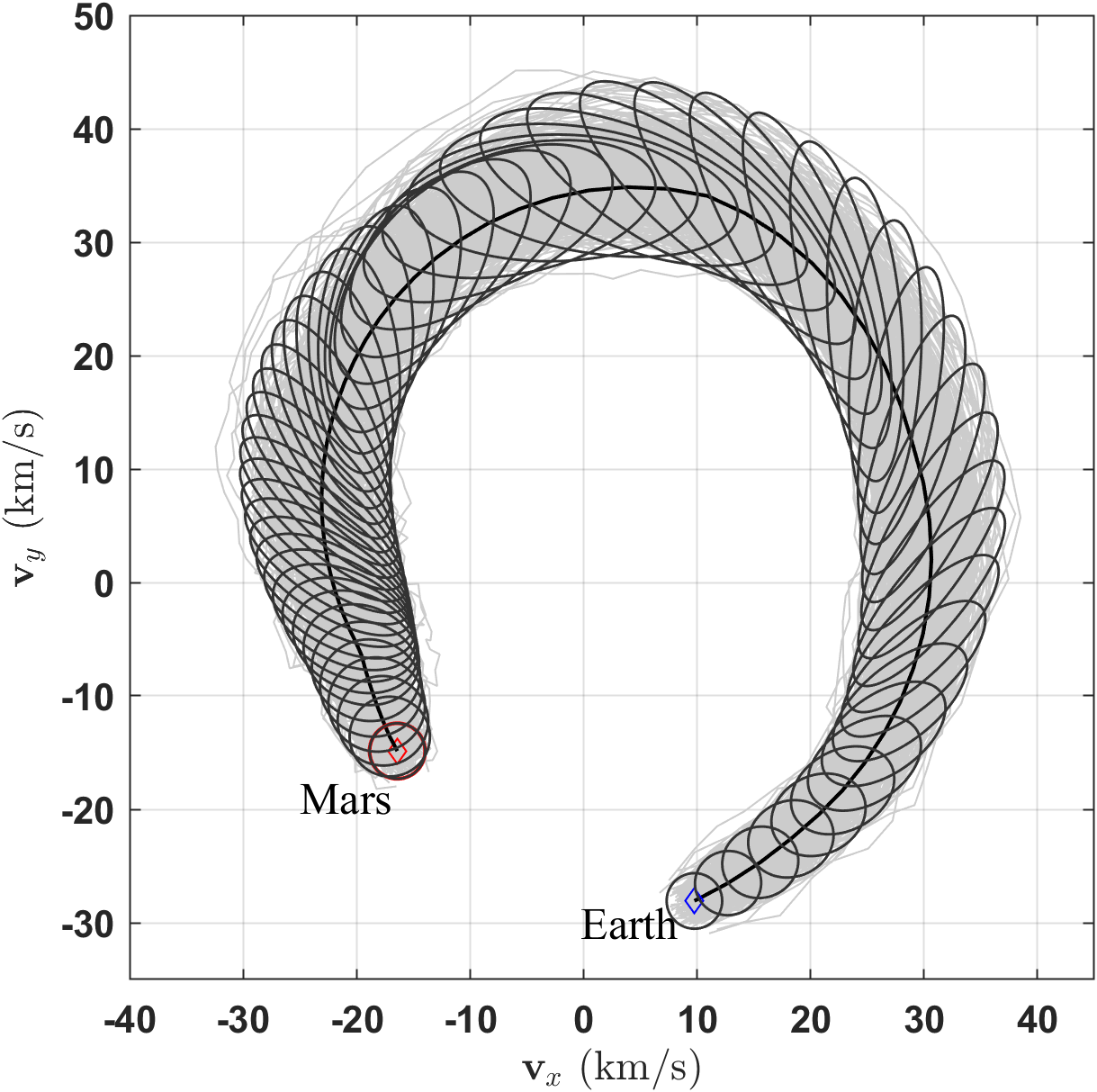}
        \caption{Closed-loop velocity trajectory and covariances (10x).}
        \label{fig:3DvelS10}
    \end{subfigure}
    \caption{3D case: Monte Carlo validation (1000 samples) within the $95\%$-confidence ellipses.}
    \label{fig:3DCLpos&VelcovS10}
\end{figure}

\begin{figure}[h!]
    \centering
    \includegraphics[width=0.45\textwidth]{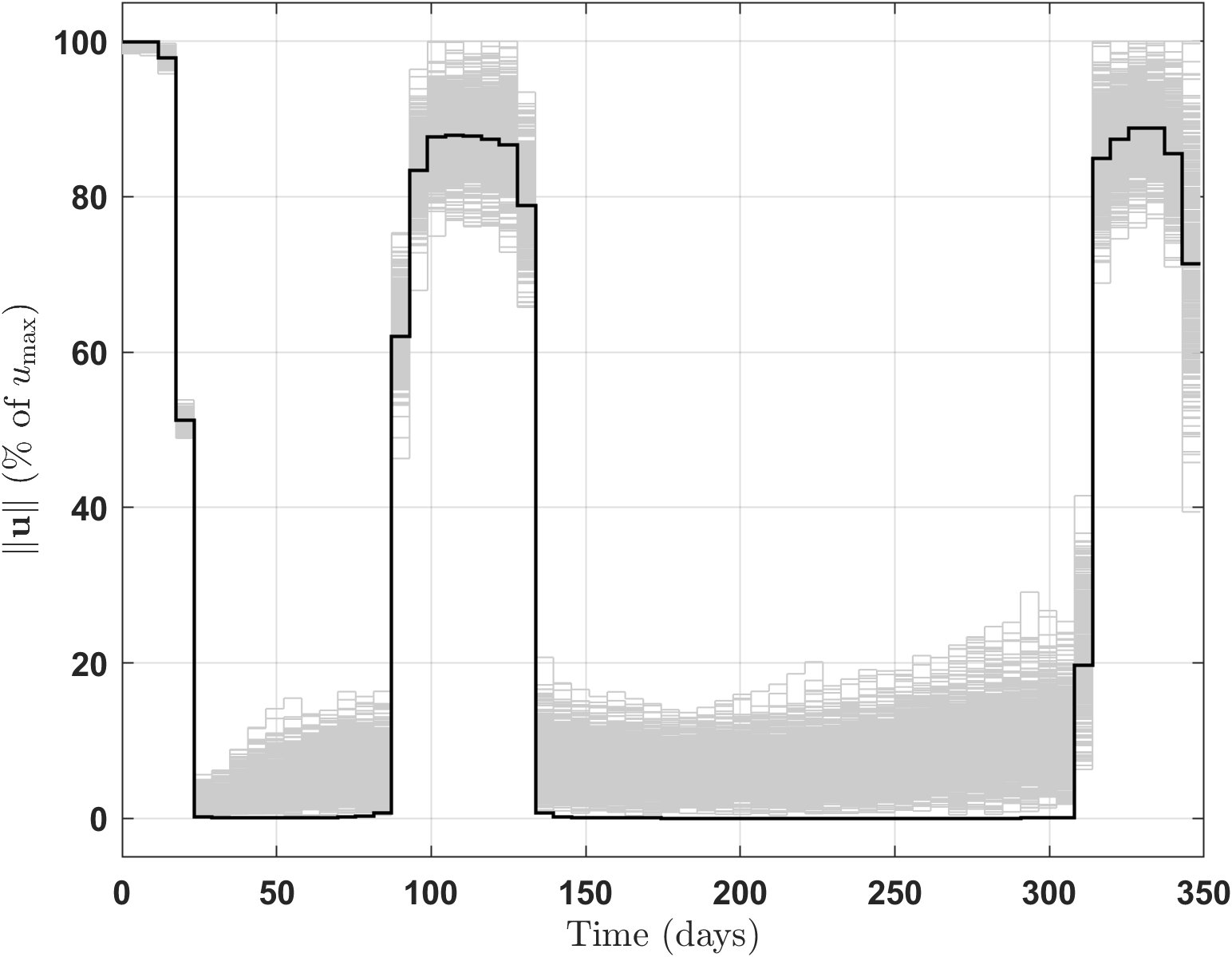}
    \caption{3D case: magnitude of the Monte Carlo control samples (gray) and feed-forward  nominal control (black).}
    \label{fig:3d_control}
\end{figure}

Figure \ref{fig:3DCLpos&VelcovS10} shows the $x$-$y$ projection of the closed-loop solution, where the position and velocity covariance ellipses are scaled by a factor of 10 for visual clarity. The close alignment of the Monte Carlo samples with the predicted $95\%$-confidence ellipses confirms the fidelity of the linearized covariance dynamics within the SCP loop. The control profile, shown in Figure \ref{fig:3d_control}, consists of three primary thrust arcs. Similar to the 2D case, the nominal thrust magnitude shows a modulated profile, deviating from a pure bang-bang structure to maintain the required robustness margins.
The control covariance ellipsoids remain within the thrust upper limit, confirming chance constraint satisfaction.

\begin{figure}[h!]
    \centering
    \includegraphics[width=0.46\textwidth]{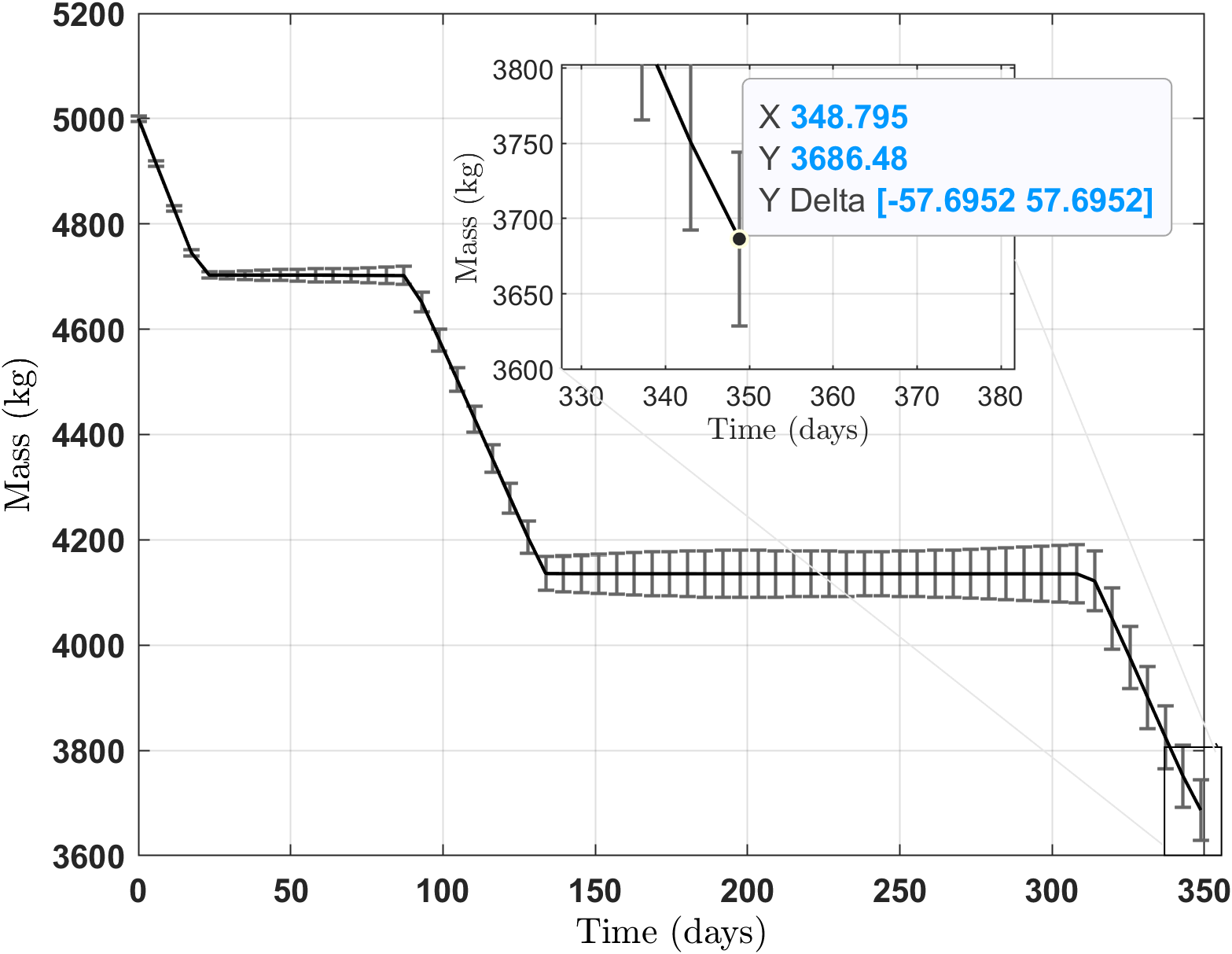}
    \caption{3D case: mass and its standard deviation vs. time. The standard deviation changes  from an initial value of zero.}
    \label{fig:mass_evolution}
\end{figure}

One of our main analyses involves a direct comparison between the full stochastic-mass model and a simplified model where mass is treated deterministically by using only the feed-forward control.
Figure \ref{fig:mass_evolution} shows the time history of the mass and its standard deviation, which is not necessarily a monotonically increasing value. In addition, the largest standard deviation occurs, in this case, at some point along the trajectory. If we enforce any constraint on the mass standard deviation at the final time to be greater than 57.7 kg ($\sigma_{m_f} \ge 57.7$ kg), this constraint is not going to become active, and the solution is more constrained by the final position and velocity variances. The mean mass decreases to a final value of approximately $m_f = 3686.48$ kg upon arrival at Mars. As expected, the final mass is less than its 2D counterpart. The mass standard deviation grows from zero to a final value of nearly 57.7 kg, driven entirely by the cumulative uncertainty in thrust execution.

\begin{figure}[h!]
    \centering
    \includegraphics[width=0.45\textwidth]{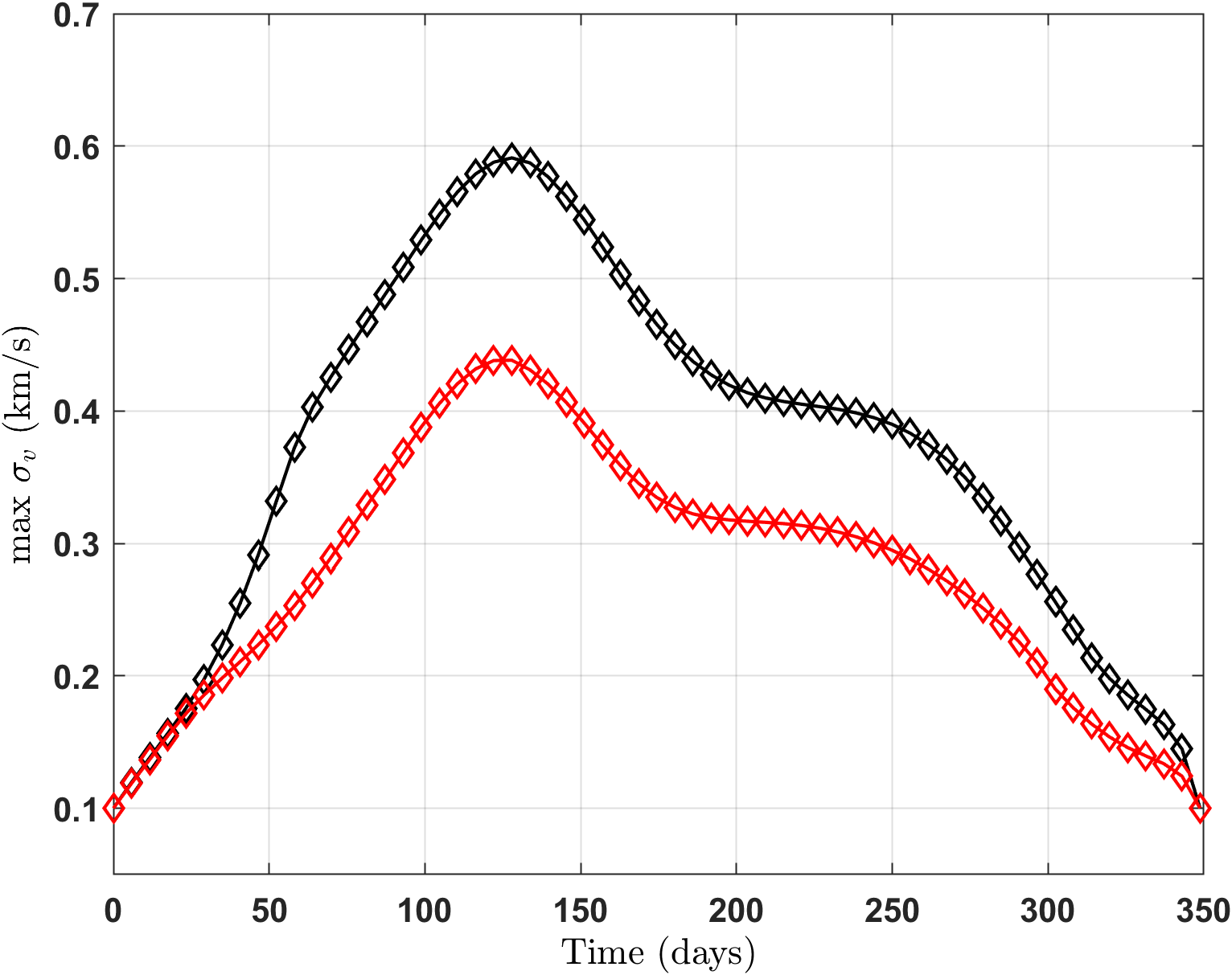}
    \caption{3D case: Comparison of velocity standard deviations. The model with mass uncertainty (black) shows significantly higher dispersion, especially during the middle thrust arc.}
    \label{fig:velocity_std_comparison}
\end{figure}

\begin{figure}[h!]
    \centering
    \includegraphics[width=0.45\textwidth]{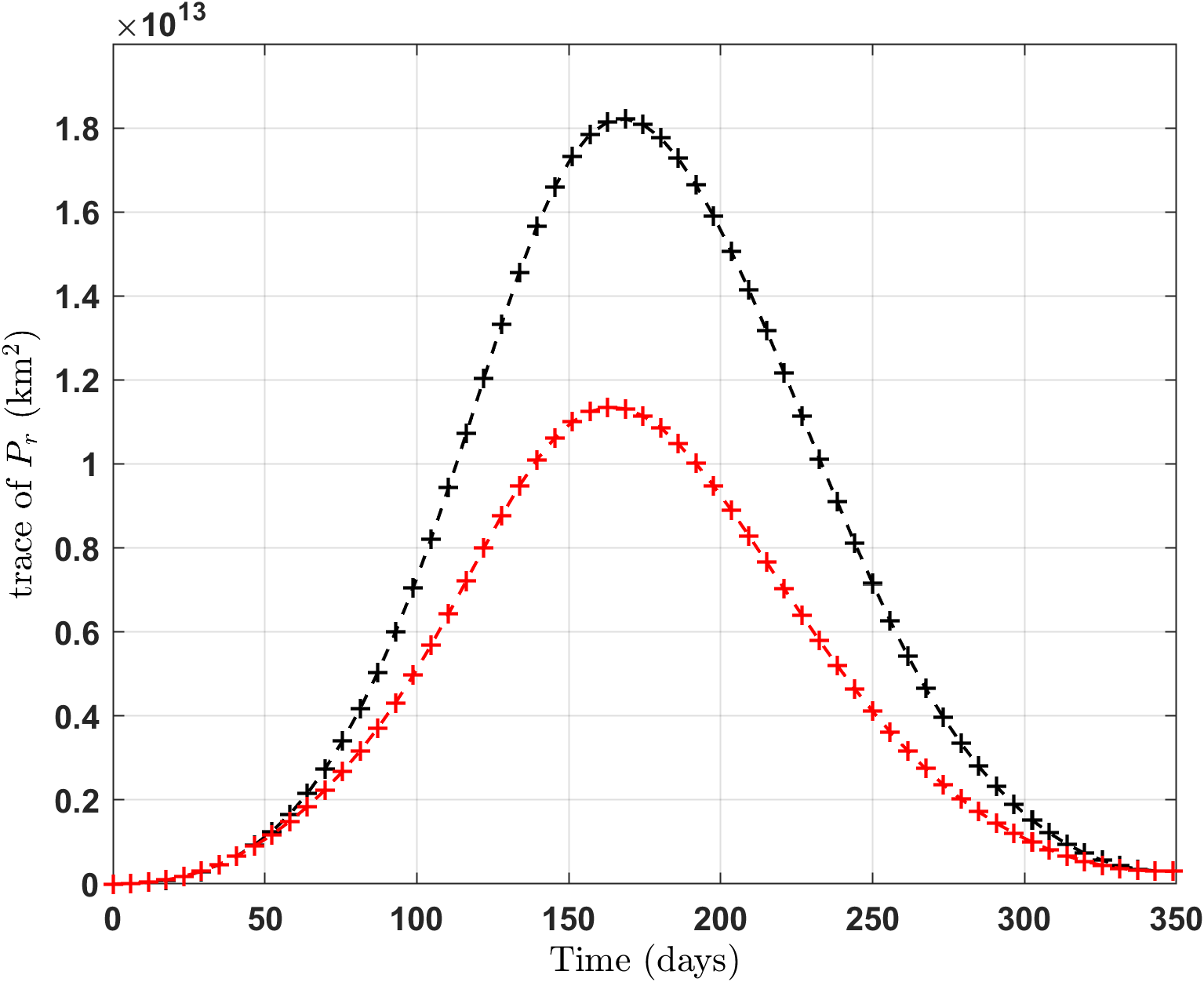}
    \caption{3D case: trace of the position covariance. The formulation with mass uncertainty (black) results in up to $61.03\%$ larger dispersion than the deterministic-mass model (red).}
    \label{fig:position_trace_comparison}
\end{figure}

This mass inclusion and its associated uncertainty has a notable impact on the translational states.
Figure \ref{fig:velocity_std_comparison} shows the velocity standard deviation along the first principal component.
The model with mass uncertainty (black) exhibits a peak velocity dispersion up to $34.84\%$ higher than the deterministic-mass model (red). This inflation occurs because mass uncertainty directly translates into acceleration uncertainty via the $1/m_t$ term in the dynamics. The effect on the overall position dispersion is summarized in Figure \ref{fig:position_trace_comparison}, which plots the trace of the position covariance matrix. Inclusion of mass uncertainty results in a $61.03\%$ higher peak position dispersion during the mid-course phase.
This notable increase highlights that neglecting mass uncertainty leads to a significant underestimation of the required control authority and the true mission risk quantified in terms of the peak position and velocity variances.

\begin{figure}[h!]
    \centering
    \includegraphics[width=0.46\textwidth]{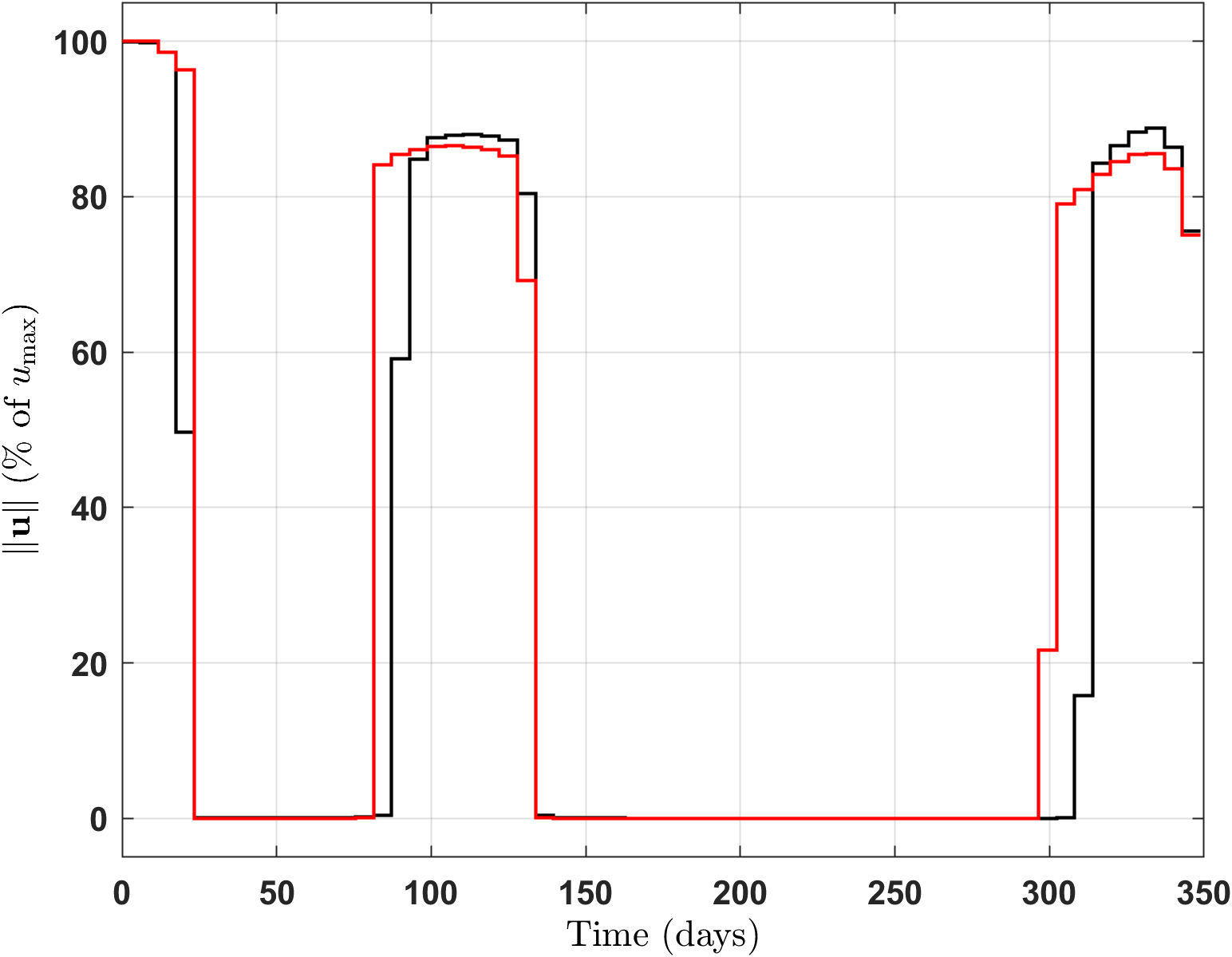}
    \caption{3D case: comparison of the feed-forward control, $\| \bm{F} \|$. The controller accounting for mass uncertainty (black) requires higher thrust to achieve the same terminal covariance bounds compared to the controller that neglects it (red).}
    \label{fig:control_comparison}
\end{figure}

Figure \ref{fig:control_comparison} compares the nominal control profiles. The controller for the full stochastic-mass model (black) commands up to $6\%$ higher thrust during the main arcs compared to the deterministic-mass controller (red). This additional control effort is necessary to counteract the larger predicted dispersions and actively shape the covariance, confirming that robustness to mass uncertainty carries a tangible fuel and performance cost. Also, the middle thrust arc leverages the Oberth effect \cite{bowerfind2024application,bowerfind2024rapid} (i.e., it utilizes the gravitational well of the Sun to maximize the rate of change of its
energy). Thus, it has a notable influence on the resulting position and velocity  and their evolution along the trajectory. 

\begin{figure}[h!]
    \centering
    \includegraphics[width=0.46\textwidth]{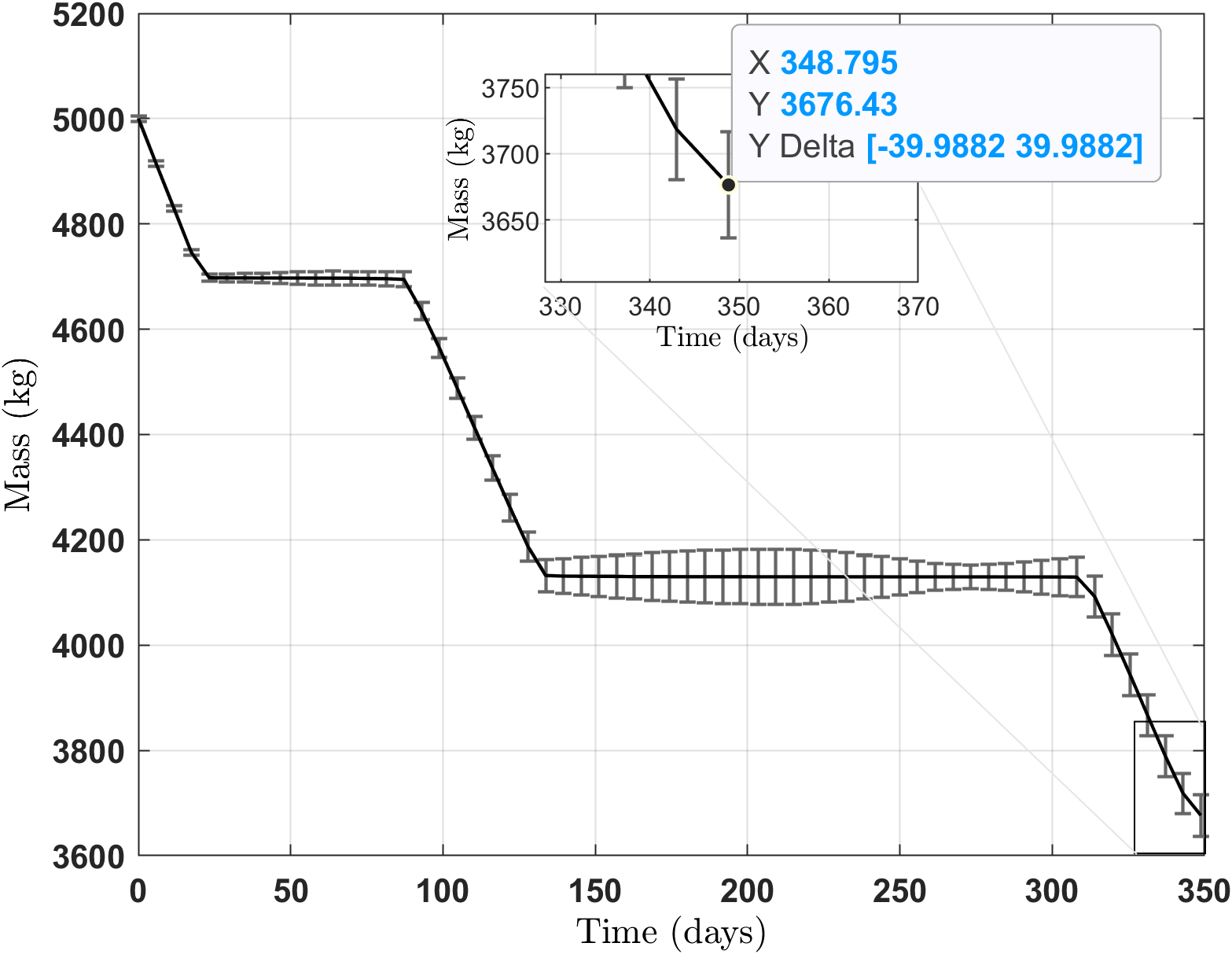}
    \caption{3D case: time history of spacecraft mass and its standard deviation with $\sigma_{m_f} \leq 40$ kg. }
    \label{fig:masssSTD40kg}
\end{figure}

\begin{figure}[h!]
    \centering
    \includegraphics[width=0.46\textwidth]{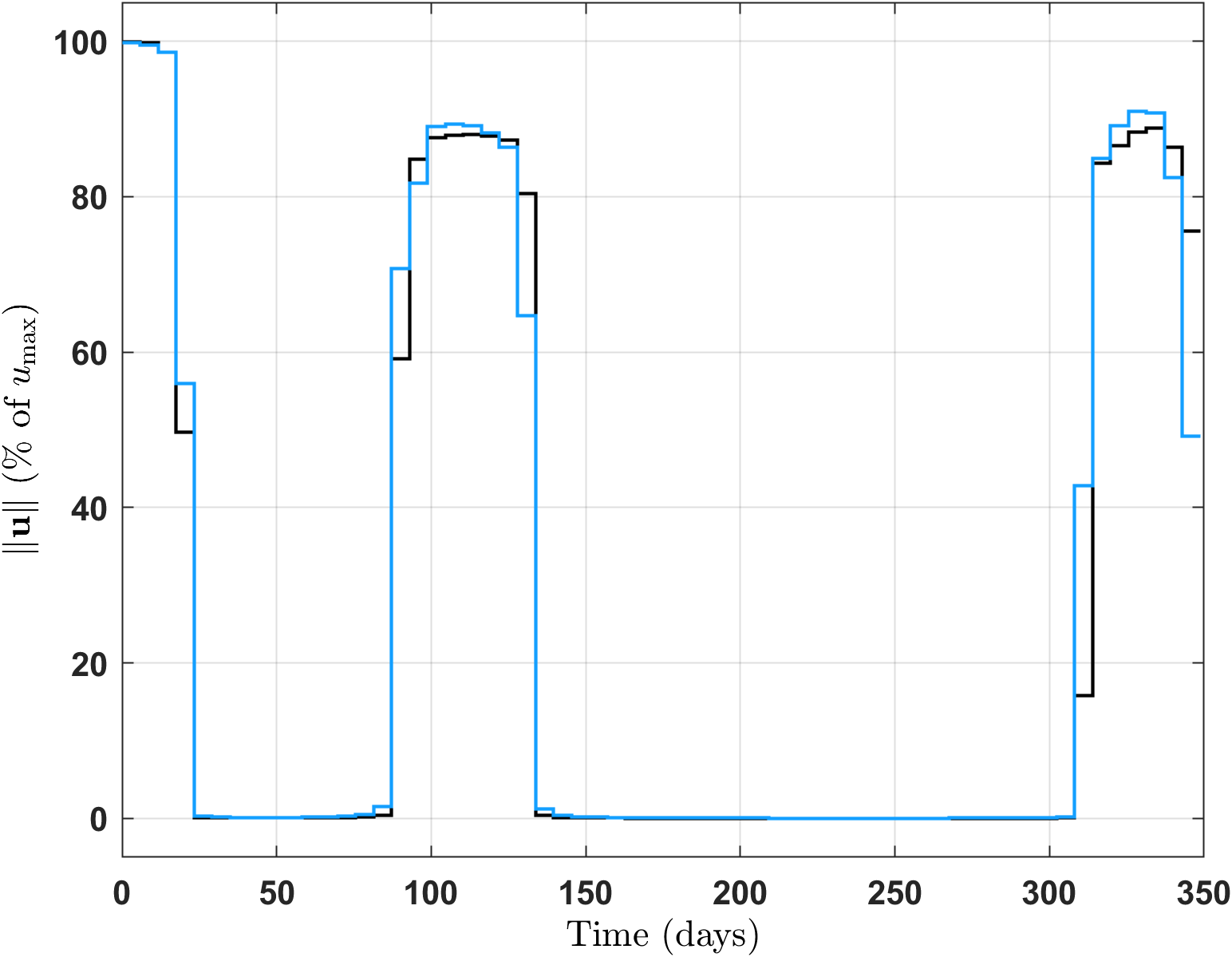}
    \caption{Comparison of nominal control thrust magnitudes.
The controller accounting for constraint with $\sigma_{m_f} \leq  40$ kg (blue) requires higher thrust levels.}
    \label{fig:CntrlMassEffect}
\end{figure}

 To further illustrate the utility of the proposed formulation, we solve the same 3D problem, but with 
\begin{align*}
P_{f} = \mathrm{diag}(10^{11},\ 10^{11},\ 10^{11},\ 10^{-2},\ \nonumber \\ 10^{-2},\ 10^{-2},\ 1600) \quad \text{(km}^2, \text{(km/s)}^2, \text{kg}^2),
\end{align*}
which corresponds to $\sigma_{m_f} \leq 40$ kg (i.e., a maximum 40 kg standard deviation in final mass) at the end of the trajectory. Figure \ref{fig:masssSTD40kg} shows the spacecraft mass and its standard deviation vs. time. An interesting observation is that while the standard deviation has decreased compared to the previous case, the final mass (3676.43) has decreased by 10 kg. More propellant is consumed to achieve the reduction of the dispersion in mass at the final time. We highlight that this result is similar to those obtained by solving optimal control problems using a desensitized optimal control formulation \cite{jawaharlal2024reduced}, which is another powerful method for generating robust optimal solutions \cite{seywald2024desensitized,chadalavada2026desensitized}. Figure \ref{fig:CntrlMassEffect} also compares the nominal thrust magnitude profiles for the previous case and when we enforce $\sigma_{m_f} \leq  40$ kg. The difference in the two profiles leads to the 10 kg more propellant consumption.

The considered cases conclusively demonstrate two key findings. First, the covariance-variable formulation is highly effective for robust trajectory optimization in complex, high-dimensional astrodynamics problems. Second, and more critically, the propagation of mass uncertainty has a substantial, quantifiable impact on both state dispersion and control effort. Neglecting this coupling, as is common in previous models in the literature, leads to an underestimation of mission risk. The proposed framework successfully manages this coupled uncertainty, providing a robust and high-fidelity solution for generating practical interplanetary guidance and control algorithms.

\section{Conclusions} \label{sec:con}
This work presents a sequential convex programming framework for solving the chance-constraint, covariance-steering problem for designing robust low-thrust interplanetary trajectories. The method incorporates a piecewise affine state-feedback policy and explicitly models spacecraft mass dynamics within the stochastic formulation. By coupling the control input and disturbance intensity through the time-varying spacecraft mass, the framework captures the evolution of uncertainty more realistically than conventional approaches. It is also possible to enforce a desired variance on the spacecraft mass at the end of the trajectory. The advantage of this model is that we can consider realistic engine parameters (such as thrust magnitude and specific impulse) and enforce constraints directly on the maximum thrust that is to be produced by the propulsion system. 

Numerical results successfully demonstrate the regulation of both the mean trajectory and the state covariance in two- and three-dimensional scenarios, producing closed-loop solutions that satisfy probabilistic terminal constraints and remain consistent with Monte Carlo dispersion analysis. In particular, the inclusion of mass as a stochastic state reveals critical differences in both control effort and uncertainty evolution that are overlooked when mass is treated deterministically.
Future studies will apply the computationally efficient proposed method to more challenging atmospheric \cite{nurre2025constrained} and cislunar \cite{saloglu2024acceleration,nurre2024end} trajectory optimization problems. 

\bibliographystyle{IEEEtran}
\bibliography{mybibfile}

\end{document}